\documentclass[a4paper]{article}

\usepackage[dvips]{graphicx}
\usepackage{dsfont}
\usepackage{amssymb}
\usepackage{amsthm}
\usepackage{latexsym}
\usepackage{amsmath}
\usepackage{color}
\usepackage{comment}
\usepackage{ascmac}
\usepackage{enumerate}
\usepackage{bm}
\usepackage{bbm} 
\usepackage{amsfonts}

\newif\ifrs
\rstrue
\ifrs \usepackage{mathrsfs} \fi

\newif\ifcol
\colfalse
\ifcol
\newcommand{\colorr}{\color[rgb]{0.8,0,0}}

\newcommand{\colorb}{\color[rgb]{0,0,0.8}}

\newcommand{\colorn}{\color[rgb]{1,1,1}}

\newcommand{\coloro}{\color[rgb]{1,0.4,0}}
\newcommand{\coloroy}{\color[rgb]{1,0.95,0}}
\newcommand{\colorsb}{\color[rgb]{0,0.95,1}}
\else
\newcommand{\colorr}{\color{black}}

\newcommand{\colorn}{\color{black}}
\newcommand{\coloro}{\color{black}}
\newcommand{\colo}{\color[rgb]{1,0.851,0}}
\newcommand{\coloroy}{\color{black}}
\newcommand{\colorsb}{\color{black}}
\fi

\newcommand{\coloraka}{\color{black}}

\excludecomment{en-text}
\includecomment{jp-text}

\excludecomment{comment}

\newtheorem{lemma}{Lemma}

\newtheorem{theorem}{Theorem}
\newtheorem{remark}{Remark}

\def\D{{\bf D}}

\def\F{{\bf F}}

\def\cale{{\cal E}}
\def\calf{{\cal F}}

\def\cali{{\cal I}}

\def\call{{\cal L}}

\def\calo{{\cal O}}

\def\cals{{\cal S}}

\def\sskip{\hspace{0.5cm}}
\def\simleq{ \raisebox{-.7ex}{\em $\stackrel{{\textstyle <}}{\sim}$} }

\def\ep{\epsilon}
\def\half{\frac{1}{2}}
\def\iku{\rightarrow}

\def\up{\uparrow}

\def\y{\vspace*{3mm}\\}

\def\nn{\nonumber}
\def\be{\begin{equation}}
\def\ee{\end{equation}}
\def\bea{\begin{eqnarray}}
\def\eea{\end{eqnarray}}
\def\beas{\begin{eqnarray*}}
\def\eeas{\end{eqnarray*}}
\def\bi{\begin{itemize}}
\def\ei{\end{itemize}}

\def\bd{\begin{description}}
\def\ed{\end{description}}
\def\l{\left}
\def\r{\right}

\newcommand{\bbA}{{\mathbb A}}

\newcommand{\bbD}{{\mathbb D}}

\newcommand{\bbH}{{\mathbb H}}
\newcommand{\bbI}{{\mathbb I}}

\newcommand{\bbN}{{\mathbb N}}

\newcommand{\bbP}{{\mathbb P}}

\newcommand{\bbR}{{\mathbb R}}

\newcommand{\bbW}{{\mathbb W}}

\newcommand{\bbZ}{{\mathbb Z}}

\def\dla{\langle\!\langle}
\def\dra{\rangle\!\rangle}
\def\koko{{\colo{koko}}}
\def\bd{\begin{description}}
\def\ed{\end{description}}

\def\D2{\bbD_{2,\infty-}}

\def\dotc{\stackrel{\circ}{C}}
\def\dotf{\stackrel{\circ}{F}}
\def\dotw{\stackrel{\circ}{W}}

\def\tj{_{t_j}}
\def\tjm{_{t_{j-1}}}
\def\intj{\int_{t_{j-1}}^{t_j}}
\def\intjm{\int_{t_{j-1}}}
\def\one{^{(1)}}
\def\two{^{(2)}}

\begin{document}

\title{Asymptotic expansion for the quadratic form of the diffusion process
\footnote{
This work was in part supported by 
Japan Society for the Promotion of Science Grants-in-Aid for Scientific Research 
No. 19340021 (Scientific Research), 
No. 24340015 (Scientific Research), 
No. 24650148 (Challenging Exploratory Research); 
the Global COE program ``The Research and Training Center for New Development in Mathematics'' 
of the Graduate School of Mathematical Sciences, University of Tokyo; 
and by a Cooperative Research Program of the Institute of Statistical Mathematics. 
}
\author{Nakahiro YOSHIDA\\
\begin{small}
University of Tokyo
\end{small}\\
\begin{small}
Graduate School of Mathematical Sciences, 
3-8-1 Komaba, Meguro-ku, Tokyo 153, Japan. 
\end{small}\\
\begin{small}
e-mail: nakahiro@ms.u-tokyo.ac.jp
\end{small}
}}
\date{
August 28, 2012
}
\maketitle
\ \\
{\it Summary} 
In \cite{Yoshida2010}, asymptotic expansion of the martingale with mixed normal limit was provided. 
The expansion formula is expressed by the adjoint of a random symbol with coefficients 
described by the Malliavin calculus, 
differently from the standard invariance principle. 
As an application, an asymptotic expansion for a quadratic form of a diffusion process was 
derived in the same paper. 
This article gives some details of the derivation,  
after a short review of the martingale expansion in mixed normal limit. 
\ \\
\ \\
{\it Keywords and phrases} 
Asymptotic expansion, martingale, mixed normal distribution, 
Malliavin calculus, random symbol, 
double It\^o integral, quadratic form. 
\ \\
\ \\
{\it Mathematics Subject Classification} 
62E20, 
60H07, 
60G44, 
62M99.  


\section{Introduction}

The quadratic form of the increments of a diffusion type process under finite time horizon, the ``realized volatility'' in financial context for example, 
is  in general asymptotically mixed normal. 
When the limit is normal, the asymptotic expansion 
of the quasi-likelihood type estimator 
was derived in \cite{Yoshida1997} as an application of the martingale expansion. 
The expansion for the quadratic form with asymptotically mixed normal limit is then indispensable 
to develop the higher-order approximation and inference for the volatility. 
However, the classical approaches in limit theorems, where the limit is a process with independent increments, 
do not work.  

The previous paper \cite{Yoshida2010}\footnote{A simplified version \cite{Yoshida2012-mart_exp} is also available now. } presented asymptotic expansion of the martingale with mixed normal limit. 
The expansion formula is expressed by the adjoint of a random symbol with coefficients 
described by the Malliavin calculus, 
differently from the standard invariance principle. 
As an application, an asymptotic expansion for a quadratic form of a diffusion process was derived in \cite{Yoshida2010}. 
The aim of this article is to give a short review of this result and some details of its derivation.

\section{Asymptotic expansion of a quadratic form of a diffusion process}

We consider a diffusion process 
satisfying the It\^o stochastic integral equation 
\bea\label{240221-1} 
X_t &=& X_0+\int_0^t b(X_s)ds+\int_0^t\sigma(X_s)dw_s . 
\eea
Here
$b$ and $\sigma$ are assumed to be smooth functions with bounded derivatives of positive order. 
We only treat one-dimensional case for notational simplicity, however, multivariate analogue is straightforward. 
Extension to general It\^o processes is also possible but the descriptions would become involved. 
We will consider a quadratic form 
\bea\label{240221-10} 
U_n &=&
\sum_{j=1}^n c(X_{t_{j-1}}) (\Delta_j X)^2, 
\eea
of the increments of $X$, 
where $t_j=j/n$ and $\Delta_jX=X_{t_j}-X_{t_{j-1}}$. 
The function $c$ is in $C^\infty_\up(\bbR)$.\footnote{$C^\infty_\up(\bbR^d;\bbR^k)$ is 
the set of $\bbR^k$-valued smooth functions defined on $\bbR^d$ 
with all derivatives of at most polynomial growth. $C^\infty_\up(\bbR;\bbR)$ is simply denoted by $C^\infty_\up(\bbR)$. }

The quadratic form (\ref{240221-10}) of the increments of $X$ appears in applications  
in statistics and finance. In the high-frequency sampling of $n$ tending to $\infty$, 
$U_n$ converges in probability to 
\beas 
U_\infty &=& \int_0^1 c(X_t)\sigma(X_t)^2 dt.
\eeas
The normalized error is 
\bea\label{240801-1}
Z_n&=& 
\sqrt{n}(U_n-U_\infty).
\eea
It is well known that $Z_n$ has a mixed normal limit distribution in general. 
However, the limit theorem is not always sufficient for approximation nor for theoretical statistics. 
Our interest is in more precise approximation to the distribution of $Z_n$.

We write $f_t$ for $f(X_t)$, given function $f$. 
For differentiable $f$, the It\^o decomposition of $f_t=f(X_t)$ is 
denoted by 
\beas 
f_t &=& f_0 + \int_0^t f^{[1]}_s dw_s
+\int_0^t f^{[0]}_sds. 
\eeas
Obviously, 
\beas 
f^{[1]}_t = \sigma(X_t)\partial_xf(X_t)\sskip{\rm and}\sskip
f^{[0]}_t = Lf(X_t)\sskip{\rm with}\sskip 
L=b\partial_x+\half\sigma^2\partial_x^2.
\eeas


For a $d_1$-dimensional reference variable, we will consider 
\bea\label{240627-1} 
F_n = 
{\coloro{
\frac{1}{n}\sum_{j=1}^n \beta(X_{t_{j-1}})
}}
\hspace{5mm}\mbox{or} \hspace{5mm}
F_n = 
F_\infty:=
{\coloro{
\int_0^1\beta(X_t) dt,
}}
\eea
where $\beta\in C^\infty_\up(\bbR;\bbR^{d_1})$. 
The results will be the same in these cases up to the first order asymptotic expansion we will discuss in this paper. 
It is standard in theoretical statistics to treat 
the joint distribution of $Z_n$ and $F_n$ because $F_n$ can be the quadratic variation of the score martingale 
and the LAMN property is then established on the joint convergence. The studentization also motivates the joint convergence. 
%

Let $a(x)=c(x)\sigma(x)^2$. Define $1+d_1$-dimensional vector fields $V_0$ and $V_1$ by  
\beas 
V_0(x_1,x_2) = 
\l[\begin{array}{c} b(x_1)-\half \sigma(x_1)\partial_{x_1}\sigma(x_1)\y
\beta(x_1) \end{array} \r]
&\mbox{  and  }&
V_1(x_1,x_2)=
\l[\begin{array}{c} \sigma(x_1)\y
0 \end{array} \r]
\eeas
for $x_1\in\bbR$ and $x_2\in\bbR^{d_1}$. 
Denote by $\mbox{Lie}[V_0;V_1](x_1,x_2)$ the Lie algebra generated by 
\beas 
V_1, \>[V_i,V_j]\>(i,j=0,1),\>[V_i,[V_j,V_k]]\>(i,j,k=0,1),....
\eeas
at $(x_1,x_2)$, where $[\cdot,\cdot]$ is the Lie bracket. 

Assume that the support $\mbox{supp}\>P^{X_0}$ of the law of $X_0$ is compact. 
We will assume the following non degeneracy conditions. 
\bd
\item[[$\bbH1$$\hspace{-2.0mm}$]]\hspace{5mm}
$\inf_{x\in\bbR} |a(x)| > 0.$
\ed
\bd
\item[[$\bbH2$$\hspace{-1.8mm}$]]\hspace{5mm}
$\mbox{Lie}[V_0;V_1](X_0,0)=\bbR^{1+d_1}$ a.s. 
\ed

The asymptotic expansion formula will be described with certain random symbols. 
The full random symbol, denoted by $\sigma(z,iu,iv)$ for $(z,u,v)\in\bbR\times\bbR\times\bbR^{d_1}$, 
consists of the adaptive random symbol $\underline{\sigma}$ and the anticipative random symbol $\overline{\sigma}$. 
These random symbols are defined as follows. 

Let
\beas 
h_t&=&
c_tb_t^2+c_tb^{[1]}_t\sigma_t-\half c^{[0]}_t\sigma_t^2-c^{[1]}_t\sigma_t\sigma^{[1]}_t
\eeas 
and 
{\coloraka 
\beas 
k_t &=& 2c_tb_t\sigma_t+c_t\sigma_t\sigma^{[1]}_t-\half c^{[1]}_t\sigma_t^2.
\eeas
}
In the present situation, the adaptive random symbol is given by 
\bea\label{240513-1}
\underline{\sigma}(z,iu,iv)
&=&
{\coloro{\frac{2z}{3}\int_0^1a(X_s)^3ds\Big( \int_0^1a(X_s)^2ds\Big)^{-1}}} \>(iu)^2 
{\coloraka +iu\int_0^1 k_t dw_t}
+iu\int_0^1 h_tdt.
\eea

The processes $D_sX_t$ and $D_rD_sX_t$ satisfy a system of partially linear equations: 
\beas 
D_sX_t &=&
\sigma(X_s)+\int_s^t {\coloraka b}'(X_t)D_sX_{t_1}d{t_1}+\int_s^t\sigma'(X_{t_1})D_sX_{t_1}dw_{t_1}
\eeas
for $s\leq t$, and  
\beas
D_rD_sX_t &=&
\sigma'(X_s)D_rX_s
+\int_s^t{\coloraka b}''(X_{t_1})D_rX_{t_1}D_sX_{t_1}d{t_1}
+\int_s^t{\coloraka b}'(X_{t_1})D_rD_sX_{t_1}d{t_1}
\\&&
+\int_s^t\sigma''(X_{t_1})D_rX_{t_1}D_sX_{t_1}dw_{t_1}
+\int_s^t\sigma'(X_{t_1})D_rD_sX_{t_1}dw_{t_1}
\eeas
for $r{\coloraka <} s\leq t$. 
The $L^p$-estimates of the solution are at hand. 
%
\begin{comment}
\beas
\sigma_{s,r}(iu,iv)
&=&
\half D_rC^{{\colorr{\infty}}}_s[u^{\otimes2}]\otimes 
\Big(iD_sW_\infty[u]-\half D_sC_\infty[u^{\otimes2}]+iD_sF_\infty[v] \Big)
\\&&
+
\Big(iD_rW_\infty[u]-\half D_rC_\infty[u^{\otimes2}]+iD_rF_\infty[v] \Big)
\otimes 
\Big(iD_sW_\infty[u]-\half D_sC_\infty[u^{\otimes2}]+iD_sF_\infty[v] \Big)
\\&&
+
\Big(iD_rD_sW_\infty[u]-\half D_rD_sC_\infty[u^{\otimes2}]+iD_rD_sF_\infty[v] 
\Big). 
\eeas
\end{comment}
%
Then the anticipative random symbol is given by the formula 
\bea\label{26720828-1}  
\bar{\sigma}(iu,iv)
&=&
 \int_0^1 
iu\>a(X_s)
\sigma_{s,s}(iu,iv)\>ds
\eea
with
\bea\label{26720828-2}  
\sigma_{s,s}(iu,iv)
&=&
\Big(-u^2\int_s^1\alpha'(X_t)D_sX_tdt+i\int_s^1\beta'(X_t)[v]D_sX_tdt \Big)^2
\nn\\&&
-u^2\int_s^1\{\alpha''(X_t)(D_sX_t)^2
+\alpha'(X_t)D_sD_sX_t\}dt
\nn\\&&
+i\int_s^1\{\beta''(X_t)[v](D_sX_t)^2
+\beta'(X_t)[v]D_sD_sX_t\}dt,
\eea
where $\alpha(x)=a(x)^2$. 
Let $C_\infty=2\int_0^1\alpha(X_t)dt$. 
%
\begin{comment}
\beas 
\bar{\sigma}(iu,iv)
&=&
\half\mbox{Tr}^*  \int_0^1 
\bar{K}^\infty(t,t)[iu]\otimes 
\sigma_{t,t}(iu,iv)\>{\coloroy{\mu (dt)}}. 
\eeas
\end{comment}

%
The density of the multi-dimensional normal distribution with mean vector ${\sf m}$ and variance matrix ${\sf C}$ 
is denoted by $\phi(z;{\sf m},{\sf C})$. 
With the full random symbol 
\bea\label{240731-4} 
\sigma(z,iu,iv) &=& \underline{\sigma}(z,iu,iv) +\overline{\sigma}(iu,iv), 
\eea
the density function $p_n(z,x)\in C^\infty(\bbR^{1+d_1})$ is defined by 
\bea\label{240731-3} 
p_n(z,x) &=& 
E\bigg[\phi(z;0,C_\infty)\delta_x(F_\infty)\bigg]
+\frac{1}{\sqrt{n}}E\bigg[\sigma(z,\partial_z,\partial_x)^*\bigg\{\phi(z;0,C_\infty)\delta_x(F_\infty)\bigg\}\bigg]. 
\eea
Here $\delta_x(F_\infty)$ is Watanabe's delta function (\cite{Watanabe1983}). 
The adjoint operation $\sigma(z,\partial_z,\partial_x)^*$ is defined by 
\beas
\sigma(z,\partial_z,\partial_x)^*
\bigg\{\phi(z;0,C_\infty)\delta_x(F_\infty)\bigg\}
&=&
\sum_j(-\partial_x)^{m_j}(-\partial_x)^{n_j}\bigg\{
c_j(z)\phi(z;0,C_\infty)\delta_x(F_\infty)\bigg\}
\eeas
for the random symbol $\sigma(z,\partial_z,\partial_x)$ having a representation
\beas 
\sigma(z,\partial_z,\partial_x)
&=&
\sum_jc_j(z)(iu)^{m_j}(iv)^{n_j}\sskip\mbox{(finite sum)}
\eeas
where $c_j$ are random functions of $z$, $m_j\in\bbZ_+^d$ (in the present case $d=1$) and 
$n_j\in\bbZ_+^{d_1}$. 
If $c_j$, $C_\infty$ and $F_\infty$ are smooth in Malliavin's sense and $F_\infty$ satisfies a suitable nondedeneracy condition,  
$C_\infty$ being nondegenerate as well, 
then this adjoint operation is well defined. 
{\coloraka These conditions are satisfied in the present situation, therefore $p_n(z,x)$ is well defined. }
See \cite{Yoshida2010} for details of random symbols and the adjoint operation. 

The following theorem gives an error bound for the approximate density $p_n(z,x)$.

\begin{theorem}\label{220829-3}
Suppose that {\rm [$\bbH1$]} and {\rm[$\bbH2$]} are satisfied. 
Then for any positive numbers $M$ and $\gamma$, 
\beas 
\sup_{f\in\cale(M,\gamma)} \bigg|
E\big[f(Z_n,F_n)\big]-\int_{\bbR^{1+d_1}}f(z,x)p_n(z,x)dzdx\bigg|
&=&
o\bigg(\frac{1}{\sqrt{n}}\bigg)
\eeas
as $n\to\infty$, 
where $\cale(M,\gamma)$ is the set of measurable functions $f:\bbR^{1+d_1}\to\bbR$ 
satisfying $|f(z,x)|\leq M(1+|z|+|x|)^\gamma$ for all $(z,x)\in\bbR\times\bbR^{d_1}$.  \\
\end{theorem}

\begin{en-text}
\begin{remark}\rm \koko
Under [H1], both $\mbox{ess.}\inf|\sigma(X_0)|>0$ and $\mbox{ess.}\inf|c(X_0)|>0$. 
Then [H2] is equivalent to the linear hull 
$\mbox{L}[\partial_{x_1}^i\beta(X_0);\>i\in\bbN]=\bbR^{d_1}$ a.s. 
It is rather simple but we prefer to keep [H2], which is suitable for more general form of $F_\infty$.  
\end{remark}
\end{en-text}
\begin{remark}\rm 
Condition [$\bbH1$] is usually ensured by the uniform ellipticity of the diffusion process $X_t$ 
and a reasonable choice of the estimator. 
In this sense, it is a natural assumption in statistical context. 
\end{remark}
\begin{en-text}
\begin{remark}\rm 
Consider a $(1+d_1)$-dimensional stochastic integral equation
\beas 
\check{X}_t &=& \check{X}_0 +\int_0^tV_0(\check{X}_s)ds
+\int_0^tV_1(\check{X}_s)\circ dw_s , \sskip  t\in[0,1]. 
\eeas
Then $\check{X}_1=(X_1,F_\infty)$, and the nondegeneracy condition entails 
the nondegeneracy of $F_\infty$ in paricular. 
\end{remark}
\end{en-text}

\begin{remark}\rm
The hybrid I method (a rough Monte-Carlo method in the {\coloraka first order asymptotic expansion} term)  
is useful in the application of the expansion formula to numerical approximation. 
Applications to volatility derivatives are in our scope. 
\end{remark}

\begin{remark}\rm 
In the present article, we have a {\it conditioning} variable as $F_n$. 
On the other hand, 
it is also possible to consider versions of our results without $F_n$ 
if we were interested in a single (not joint) expansion. 
It will reduce the  
differentiability conditions of variables. 
However, considering the joint distribution is natural in non-ergodic statistics. 
Studentization is important in any case. 
\end{remark}

\begin{remark}\rm 
It is also possible to obtain asymptotic expansion of the conditional distribution. 
\end{remark}

{\coloraka Section \ref{241224-1} will give some details of derivation of Theorem \ref{220829-3}. }


\section{Review of the asymptotic expansion of a double stochastic integral having a mixed normal limit}
In this section, we will give a short review of the martingale expansion. 
We refer the reader to \cite{Yoshida2010} for details. 

On a stochastic basis 
$(\Omega,\calf,\F=(\calf_t)_{t\in[0,1]},P)$ 
with $\calf=\calf_1$, 
we shall consider a sequence of $d$-dimensional functionals 
with decomposition
\bea\label{240731-1}
Z_n &=& M_n + W_n+r_n N_n. 
\eea
Here, for every $n\in\bbN$, $M^n=(M^n_t)_{t\in[0,1]}$ denotes a $d$-dimensional {\coloraka continuous} martingale 
with respect to $\F$ and $M_n = M^n_{1}$. 
In this decomposition, we assume 
$W_n$, $N_n\in\calf(\Omega;\bbR^d)$\footnote{The set of $d$-dimensional 
measurable mappings.} 
and $(r_n)_{n\bbN}$ is a sequence of positive numbers tending to zero 
as $n\to\infty$.

In this section, the reference variables are $F_n\in\calf(\Omega;\bbR^{d_1})$ ($n\in\bbN$) are general and we do not assume a specific structure like (\ref{240627-1}). 
It is possible to give asymptotic expansion of $\call\{Z_n\}$ under certain conditions; see \cite{Yoshida2010}. 
The same paper applied the expansion to the case where 
$M_n$ is given as a sum of double It\^o integrals, as reviewed in what follows. 
Let $(\bbW,\bbP)$ be an ${\sf r}$-dimensional Wiener space over time interval $[0,1]$ 
and let $H$ be the Cameron-Martin subspace of $\bbW$. 
Suppose that the probability space $(\Omega,\calf,P)$ is such that 
$\Omega=\Omega'\times\bbW$, $\calf=\calf'\times{\bf B}(\bbW)$ and $P=P'\times\bbP$ for some 
probability space $(\Omega',\calf',P')$. We will use the partial Malliavin calculus 
on $\Omega$ based on the shifts in the direction of $H$. 
For a Hilbert space $E$, 
the Sobolev space of $E$-valued functionals on $\Omega$ with indices $s\in\bbR$ 
for differentiability and $p\in(1,\infty)$ for integrability is denoted by 
$\bbD_{s,p}(E)$.

Let $\check{d}=d+d_1$ and $\ell=\check{d}+6$. 
Let 
$\dot{K}^n\in\bbD_{\ell+1,\infty}(H\otimes\bbR^d)$ and 
$\ddot{K}^n\in\bbD_{\ell+1,\infty}(H\otimes\bbR^{{\sf r}}\otimes\bbR^d)$. 
Since $H$ can be identified with $L^2([0,1];\bbR^{{\sf r}})$, 
the functionals $\dot{K}^n$ and $\ddot{K}^n$ are respectively regarded as 
$\bbR^d\times\bbR^{\sf r}$-valued and $\bbR^d\otimes\bbR^{\sf r}\otimes\bbR^{\sf r}$-valued $L^2$-functions on $[0,1]$: 
\beas 
\dot{K}^n=\big(\dot{K}^n(s)^i_\alpha\big)_{0\leq s\leq 1,\> i=1,...,d,\>\alpha=1,...,{\sf r}}
&\mbox{and}&
\ddot{K}^n=\big(\ddot{K}^n(s)^i_{\alpha\beta}\big)_{0\leq s\leq1,\>i=1,...,d,\>\alpha,\beta=1,...,{\sf r}}. 
\eeas

The sequence 
$\{t_j\}_{j=0,1,...,\bar{j}^n}$ $(n\in\bbN)$ {\coloro{is}} a triangular array 
of numbers such that $t_j=t^n_j$ depending on $n$ and that 
$0=t_0<t_1<\cdots<t_{\bar{j}^n}=1$. 
Functional $K^n(s,r)$ is defined as the $\bbR^d\otimes\bbR^{\sf r}\otimes\bbR^{\sf r}$-valued function
\beas 
K^n(s,r)
&=& 
\bigg[
r_n^{-1}\sum_j1_{(t_{j-1},t_j]}(s)\dot{K}^n(s)^i_\alpha
1_{(t_{j-1},s]}(r)\ddot{K}^n(r)^i_{\alpha\beta}\bigg]_{i=1,...,d\atop\alpha,\beta=1,...,{\sf r}}
\eeas
for $(s,r)\in[0,1]^2$. 
We write 
\beas 
\bar{K}^n(s,r) =
\big[\bar{K}^n(s,r)^i_{\alpha\beta}\big]_{i=1,...,d\atop\alpha,\beta=1,...,{\sf r}}
=
 \bigg[\dot{K}^n(s)^i_\alpha\ddot{K}^n(r)^i_{\alpha,\beta}\bigg]_{i=1,...,d\atop\alpha,\beta=1,...,{\sf r}}.
\eeas
Suppose that $\dot{K}^n$ and $\ddot{K}^n$ are progressively measurable. 
More strongly, 
{\coloro{we assume the strong predictability 
condition that $\dot{K}^n(s)$ is $\calf_{t_{j-1}}$-measurable 
for $s\in(t_{j-1},t_j]$.
Corresponding to the kernel $K^n$, we consider 
$M^n_t$ given by 
\bea\label{240731-2} 
M^n_t 
&=&
\bigg[
r_n^{-1}\sum_j \sum_{\alpha}
\int_{t_{j-1}\wedge t}^{t_j\wedge t}
\dot{K}^n(s)^i_\alpha
\Bigl(\sum_\beta\int_{t_{j-1}}^s \ddot{K}^n(r)^i_{\alpha\beta}dw_r^\beta\Bigr) dw^\alpha_s
\bigg]_{i=1,...,d}.
\eea
\begin{en-text}
We write 
\beas 
\bar{K}^n(s,r)=
\big[\bar{K}^n(s,r)^i_{\alpha\beta}\big]_{i=1,...,d\atop\alpha,\beta=1,...,{\sf r}}
=
\bigg(\dot{K}^n(s)^i_\alpha\ddot{K}^n(r)^i_{\alpha\beta}\bigg)_{i=1,...,d\atop \alpha,\beta=1,...,{\sf r}}. 
\eeas

We will consider the kernel function
%
$K^n(s,r)$ defined by 
\bea\label{211003-1}
K^n(s,r) 
&=& 
r_n^{-1}\sum_j1_{(t_{j-1},t_j]}(s)\dot{K}^n(s)
\oslash
1_{(t_{j-1},s]}(r)\ddot{K}^n(r)
\eea
Here $\oslash$ denotes the tensor product.  
That is, 
\end{en-text}

Write $C^n_t = \langle M^n \rangle_t$ and $C_n = \langle M^n \rangle_1$. 
Suppose 
$\max_j|I_j|=o(r_n)$, where 
$I_j=(t_{j-1},t_j]$, and that the sequence of measures 
\bea\label{220131-3}
\mu^n=r_n^{-2}\sum_j |I_j|^2 \delta_{t_{j-1}} 
&\to&
\mu
\eea
weakly for some measure $\mu$ on $[0,1]$ {\coloraka with a bounded derivative}. 
{\colorr We will assume that $r_n^{-8}\sum_j|I_j|^5=O(1)$. }
\begin{en-text}
\footnote{Obviously this condition is satisfied if one replaces 
$\dot{K}^n(s)$ by $\dot{K}^n(t_{j-1})$ 
in the representation of the kernel funciton $K^n(s,r)$. 
In examples, 
it will turn out that this replacement does not 
cause any practical difficulty. 
}
\end{en-text}
{\coloro{
\begin{en-text}
then 
\beas 
K^n(s,r) 
&=&
r_n^{-1}\sum_j1_{(t_{j-1},t_j]}(s)
1_{(t_{j-1},s]}(r)\bar{K}^n(s,r).
\eeas
\end{en-text}
%
\begin{en-text} 
By definition, 
the stochastic integral of $K^n(s,\cdot)$ has 
a usual meaning and 
\beas 
\int_0^1 K^n(s,r)dw_r
&=&
r_n^{-1}\sum_j1_{(t_{j-1},t_j]}(s)
\dot{K}^n(s)
\oslash\int_{t_{j-1}}^s\ddot{K}^n(r)dw_r, 
\eeas
and
\end{en-text} 
%
In this case, 
\begin{en-text}
\beas 
C_n &=& 
{\colorr{\mbox{Tr}^o}}
\int_0^1 \Bigl(\int_0^1 K^n(s,r)dw_r\Bigr)^{\otimes2} ds 
\eeas
\end{en-text}
\begin{en-text}
\beas 
C_n &=& 
{\colorr{\mbox{Tr}^o}}
\int_0^1 \bigg(\dot{K}^n(s)
\oslash
\Bigl(\int_{t_{j-1}}^s \ddot{K}^n(r)dw_r\Bigr) dw_s
\bigg)^{\otimes2} ds 
\eeas
\end{en-text}
\beas 
C_n &=& 
\bigg[
\sum_\alpha\int_0^1
r_n^{-2}\sum_j1_{I_j}(s)
\dot{K}^n(s)^{i_1}_\alpha\Bigl(\sum_\beta\int_{t_{j-1}}^s \ddot{K}^n(r)^{i_1}_{\alpha\beta}dw_r^\beta\Bigr)
\\&&\cdot 
\dot{K}^n(s)^{i_2}_\alpha\Bigl(\sum_\beta\int_{t_{j-1}}^s \ddot{K}^n(r)^{i_2}_{\alpha\beta}dw_r^\beta\Bigr)
\>ds
\bigg]_{i_1,i_2=1,...,d}
\eeas
\begin{en-text}
\beas 
C_n &=& 
\mbox{Tr}^o\int_0^1
r_n^{-2}\sum_j1_{I_j}(s)\Big(\dot{K}^n(s)\oslash
\int_{t_{j-1}}^s\ddot{K}^n(r)dw_r\Big)^{\otimes2}\>ds
\eeas
\end{en-text}
and the in-p limit of $C_n$ will be 
\begin{en-text}
\bea\label{210930-1} 
C_\infty &=& 
{\colorr{\half\>}}
\mbox{Tr}^*
\int_0^1 \bar{K}^\infty(t,t)^{\otimes2} \mu(dt)
\eea 
\end{en-text}
\bea\label{210930-1} 
C_\infty &=& 
\bigg[\>
{\colorr{\half\>}}
\sum_{\alpha,\beta}
\int_0^1 \bar{K}^\infty(t,t)^{i_1}_{\alpha\beta}\bar{K}^\infty(t,t)^{i_2}_{\alpha\beta}\>\mu(dt)
\bigg]_{i_1,i_2=1,...d}
\nn\\&=&
\half\>\mbox{Tr}^*\int_0^1 \bar{K}^\infty(t,t)\otimes\bar{K}^\infty(t,t)\>\mu(dt)
\eea 
under the conditions we will assume, where $\bar{K}^\infty$ is the limit of $\bar{K}^n$, and 
$\mbox{Tr}^*$ is the trace on $(\bbR^{\sf r}\otimes\bbR^{\sf r})^*\otimes(\bbR^{\sf r}\otimes\bbR^{\sf r})$.  
\begin{en-text}
where 
$\mbox{Tr}^o$ and 
$\mbox{Tr}^*$ denote the traces of 
the element of 
$\mbox{L}(\bbR^{\sf r};\bbR^{\sf r})$ and 
$\mbox{L}(\bbR^{\sf r}\otimes\bbR^{\sf r}
;\bbR^{\sf r}\otimes\bbR^{\sf r})$, 
respectively. 
\end{en-text}

Let $q\in(1/3,1/2)$. 
Let $\Delta=\{(s,r);\> 0\leq r\leq s\leq 1\}$  
and $\Delta^n=\cup_j\{(s,r);\>t_{j-1}\leq s \leq t \leq t_j\}$. 
{\coloro{
Let $\bbI_s=-\half(C^\infty_s-C^\infty_1)$. 

\begin{en-text}
\bd
\item[{[A1]}$^\flat$] 
Conditions in [A1] except {\rm (iii)} hold. 
\ed
\end{en-text}

\bd
\item[{[$\bbA1$]}]
{\bf (i)} 
$\bar{K}^n\in\bbD_{\ell+1}(H\otimes H\otimes \bbR^d)$ and 
a representation density of each derivative admits 
\beas 
ess.\sup_{r_1,...,r_k,\in(0,1),\atop(s,r)\in\Delta^n,\>n\in\bbN}
\bigg\|D_{r_1,...,r_k}\bar{K}^n(s,r)\bigg\|_p<\infty
\eeas
for every $p\in(1,\infty)$ and $k\leq\ell+1$. 
\begin{description}
\item[(ii)] For every $\eta>0$ and $p\in(1,\infty)$,
\beas 
\sup_{s\in(0,1),e\in S^{d-1}}
\bigg\|\bigg[\frac{\bbI_s[e^{\otimes2}]}{(1-s)^{1+\eta}}\bigg]^{-1}\bigg\|_p<\infty. 
\eeas 

\item[(iii)]\footnote{Condition [A1] (iv) of \cite{Yoshida2010}.} 
%
For every $p>1$, 
\beas
\sup_j\sup_{(s,r)\in\Delta^n\atop s\in I_j} \|\bar{K}^n(s,r)-\bar{K}^n(t_{j-1},t_{j-1})\|_{\ell,p}
&=& O(r_n^{2q})
\eeas
and
\beas 
\sup_j\sup_{t\in I_j}
\|\bar{K}^n(t_{j-1},t_{j-1})-\bar{K}^\infty(t,t)\|_{\ell,p}&=& O(r_n^{2q})
\eeas
as $n\to\infty$. 
\ed
\ed

\begin{remark}\rm 
Condition [$\bbA1$](iii) is [A1](iv) of \cite{Yoshida2010}.
In \cite{Yoshida2010}, 
[A1]$^\flat$  was ``Conditions in [A1] except {\rm (iii)} hold.'' 
\end{remark}

\begin{remark}\rm 
{\bf (i)}  
In typical cases $r_n=n^{-1/2}$ so that $r_n^{2q}=n^{-q}>n^{-\half}$ for $n>1$. 
\bd
\item[(ii)] 
Under [$\bbA1$] (i), 
for every $p\in(1,\infty)$ and $k\leq\ell$, 
\beas 
ess.\sup_{r_1,...,r_k,s\in(0,1)}
\bigg\|\frac{|D_{r_1,...,r_k}\bbI_s|}{1-s}\bigg\|_p<\infty.
\eeas

\item[(iii)] 
As for [$\bbA1$] (ii), we need 
{\coloro{the nondegeneracy of the derivative of $C^\infty_s$ in $s$, or}} 
a large deviation argument, 
in order to control $\exp(\half(C^\infty_s-C^\infty_1)[u^{\otimes2}]))$ 
for $s$ near $1$. 

\ed
\end{remark}

{\colorsb{

Let $C_\infty\in\calf(\Omega;\bbR^d\otimes\bbR^d)$, 
$W_\infty\in\calf(\Omega;\bbR^d)$ 
and 
$F_\infty\in\calf(\Omega;\bbR^{d_1})$. 
Set 
$
\dotc_n= r_n^{-1}(C_n-C_\infty)$, 
$\dotw_n= r_n^{-1}(W_n-W_\infty)$  and 
$\dotf_n = r_n^{-1}(F_n-F_\infty)$.
Consider an extention 
\beas 
(\bar{\Omega},\bar{\calf},\bar{P})=
(\Omega\times{\stackrel{\circ}{\Omega}},\calf\times{\stackrel{\circ}{\calf}},P\times{\stackrel{\circ}{P}})
\eeas
of $(\Omega,\calf,P)$ by a probability space 
$({\stackrel{\circ}{\Omega}},{\stackrel{\circ}{\calf}},{\stackrel{\circ}{P}})$. 
\begin{comment}
{\coloro{
\ \\
-------------------------
?so that 
$\zeta$ defined on $\bar{\Omega}$ depends only on $\omega'\in\Omega'$. 

(or $(M^\infty_\cdot,
N_\infty,\dotc_\infty,\dotw_\infty,\dotf_\infty,
C^\infty_\cdot,W_\infty,F_\infty)$) 
\\ 
------------------------
\ \\
}}
\end{comment}
%
Suppose that 
$M^\infty\in\calf(\bar{\Omega};C([0,1];\bbR^d))$, 
$N_\infty\in\calf(\bar{\Omega};\bbR^d)$, 
$\dotc_\infty\in\calf(\bar{\Omega};\bbR^d\otimes\bbR^d)$, 
$\dotw_\infty\in\calf(\bar{\Omega};\bbR^d)$ 
and 
$\dotf_\infty\in\calf(\bar{\Omega};\bbR^{d_1})$.

For $\check{\calf}=\calf\vee\sigma[M^\infty_1]$,  
there exists a measurable mapping 
$\check{C}_\infty:\Omega\times\bbR^d\iku\bbR^d\otimes\bbR^d$ 
such that 
\beas 
\check{C}_\infty(\omega,M_\infty)&=&E[\dotc_\infty|\check{\calf}].
\eeas
%
%
Similarly we define 
$\check{W}_\infty(\omega,z)$, $\check{F}_\infty(\omega,z)$ and 
$\check{N}_\infty(\omega,z)$ by 
\beas 
\check{W}_\infty(\omega,M_\infty)&=&E[\dotw_\infty|\check{\calf}]
\\
\check{F}_\infty(\omega,M_\infty)
&=&E[\dotf_\infty|\check{\calf}]
\\
\check{N}_\infty(\omega,M_\infty)
&=&E[N_\infty|\check{\calf}].
\eeas
Further, we introduce the notation 
\beas 
\tilde{C}_\infty(z) \equiv  
\tilde{C}_\infty(\omega,z) &:=& \check{C}_\infty(\omega,z-W_\infty)
\eeas
and similarly 
$\tilde{W}_\infty(\omega,z)$, $\tilde{F}_\infty(\omega,z)$ and 
$\tilde{N}_\infty(\omega,z)$.  
Let 
$s_n$ be a positive functional defined on $\Omega$. 
It is said that a functional $c:\Omega\times\bbR^d\to\bbR$ is a $\bbD_{s,\infty}$-polynomial 
if $c$ is a polynomial of $z\in\bbR^d$ with coefficients in $\bbD_{s,\infty}$.

\bd
\item[{[$\bbA2$]}] 
{\bf (i)}  
$F_\infty\in\bbD_{\ell+1,\infty}(\bbR^{d_1})$ 
and 
$W_\infty\in\bbD_{{\coloro{\ell+1}},\infty}(\bbR^d)$
. 

\bd
\item[(ii)]
$F_n\in\bbD_{\ell+1,\infty}(\bbR^{d_1})$, 
$W_n\in\bbD_{{\coloro{\ell+1}},\infty}(\bbR^d)$, 
$N_n\in\bbD_{{\coloro{\ell+1}},\infty}(\bbR^d)$  
and $s_n\in\bbD_{\ell,\infty}(\bbR)$. 
Moreover, 
\beas 
\sup\Big\{ 
\|\dotw_n\|_{{\coloro{\ell+1}},p}
+\|\dotf_n\|_{\ell+1,p}+\|N_n\|_{{\coloro{\ell+1}},p}+\|s_n\|_{\ell,p} 
\Big\} <\infty.
\eeas
for every $p\geq2$. 

\begin{en-text}
\item[(iii)]
{\coloro{$\lim_{n\to\infty}P\l[|\xi_n|\leq\half\r]=1$. }}

\item[(iv)]
$|C_n-C_\infty|>r_n^{1-a}$ implies $|\xi_n|\geq 1$, 
where $a\in(0,1/3)$ is a constant.

\item[(v)] For every $p\geq2$, 
\beas 
\limsup_{n\to\infty}E\Big[1_{\{|\xi_n|\leq1\}} \Delta_{(M_n+W_\infty,F_\infty)}^{-p}\Big] 
<\infty
\eeas
and moreover 
$\Delta_{F_\infty}^{-1}, \>C_\infty^{-1}\in\cap_{p\geq2} L^p$. 
\end{en-text}

\item[(iii)] 
$\tilde{C}_\infty^{j,k}$, $\tilde{W}_\infty^j$, $\tilde{N}_\infty^j$ $(j,k=1,...,d)$ 
and $\tilde{F}_\infty^l$ $(l=1,...,d_1)$ are 
$\bbD_{\ell_0,\infty}$-polynomials in $z\in\bbR^d$, where $\ell_0=2[(d_1+3)/2]$. 

\item[(iv)] 
$({\coloro{M^n_\cdot}},N_n,\dotc_n,\dotw_n,\dotf_n)
\iku^{d_s(\calf)}
(M^\infty_\cdot,N_\infty,\dotc_\infty,\dotw_\infty,\dotf_\infty)$. 

\item[(v)] 
For $G=W_\infty$ and $F_\infty$, 
\beas 
ess. \sup_{r_1,...,r_k\in(0,1)} 
\|D_{r_1,...,r_k}G\|_p 
&<&\infty
\eeas
for every $p\in[2,\infty)$ and $k\leq {\coloraka \ell+1}$. 
Moreover, $r\mapsto D_rG$ and $(r,s)\mapsto D_{r,s}G$ are 
continuous a.e. with respect to the Lebesgue measures. 
\ed
\ed

\ \\

The nondegeneracy of $(M^n_t+W_\infty,F_\infty)$ will be necessary. 
\bd
\item[{[$\bbA3$]}] 
{\bf (i)} 
There exist a sequence $(t_n)_{n\in\bbN}$ in {\colorr{$[0,1]$ with $\sup_nt_n<1$}}  
such that 
%
%

$
\sup_{t\geq t_n}
P\big[
\det\sigma_{(M^n_t+W_\infty,F_\infty)}< s_n\big] =O(r_n^\nu)
$
as $n\to\infty$ for some $\nu>\ell/3$. 

\begin{description}
\item[(ii)] 
For every $p\geq2$, 
$
\limsup_{n\to\infty} E[s_n^{-p}]<\infty  
$. 



\end{description}
\ed

\begin{remark}\rm 
The nondegeneracy $\det C_\infty^{-1}\in\cap_{p\geq2} L^p$ follows from [$\bbA1$] (ii). 
Indeed, it implies 
\beas 
\sup_{e\in\cals^{d-1}}P[C_\infty[e^{\otimes2}]<\ep] &\leq& \ep^p\sup_{e\in\cals^{d-1}}\|\bbI_0[e^{\otimes2}]^{-1}\|_p^{{\coloraka p}}
\>\leq\>C_p\ep^p\sskip(\ep>0)
\eeas
for some constant $C_p$ for every $p>0$. 
Then the desired inequality is obtained; see e.g. Lemma 2.3.1 of Nualart \cite{Nualart2006}. 
\end{remark}

\vspace{3mm}

Now we recall the martingale expansion (\cite{Yoshida2008,Yoshida2010}). 
For $Z_n$ given by (\ref{240731-1}) and (\ref{240731-2}), 
the random symbols are specified as follows. 
The adaptive random symbol is  
\bea\label{240218-1} 
\underline{\sigma}(z,iu,iv)
&=&
\half\tilde{C}_\infty(z)[(iu)^{\otimes2}]
+
\tilde{W}_\infty(z)[iu]
+
\tilde{N}_\infty(z)[iu]
+
\tilde{F}_\infty(z)[iv]. 
\eea
For the double stochastic integral in question, the anticipative random symbol is given by 
\beas 
\bar{\sigma}(iu,iv)
&=&
\half\mbox{Tr}^*  \int_0^1 
\bar{K}^\infty(t,t)[iu]\otimes 
\sigma_{t,t}(iu,iv)\>{\coloroy{\mu (dt)}},
\eeas
where the random symbol $\sigma_{t,t,}(iu,iv)$ has the expression 
\beas 
\sigma_{t,t}(iu,iv)
&=&
\Big(iD_tW_\infty[u]-\half D_tC_\infty[u^{\otimes2}]+iD_tF_\infty[v] \Big)
\otimes 
\Big(iD_tW_\infty[u]-\half D_tC_\infty[u^{\otimes2}]+iD_tF_\infty[v] \Big)
\\&&
+
\Big(iD_tD_tW_\infty[u]-\half D_tD_tC_\infty[u^{\otimes2}]+iD_tD_tF_\infty[v] 
\Big) 
\eeas
with the representation densities of the Malliavin derivatives of functionals. 
The approximate density $p_n(z,x)$ is defined by 
\bea\label{240731-5} 
p_n(z,x) &=& 
E\bigg[\phi(z;W_\infty,C_\infty)\delta_x(F_\infty)\bigg]
+r_nE\bigg[\sigma(z,\partial_z,\partial_x)^*\bigg\{\phi(z;W_\infty,C_\infty)\delta_x(F_\infty)\bigg\}\bigg]. 
\eea
for the full random symbol 
\beas 
\sigma=\underline{\sigma}+\overline{\sigma}.  
\eeas
Note that $u$ is $d$-dimensional here. 
Class $\cale(M,\gamma)$ will be abused for functions on $\bbR^{\check{d}}$.  

\begin{theorem}\label{26720826-2}
Suppose that Conditions {\rm [$\bbA1$]}, {\rm[$\bbA2$]} and {\rm [$\bbA3$]} are fulfilled. 
Then for any positive numbers $M$ and $\gamma$, 
\beas
\sup_{f\in\cale(M,\gamma)} \bigg|
E\big[f(Z_n,F_n)\big]-\int_{\bbR^{\check{d}}}f(z,x)p_n(z,x)dzdx\bigg|
&=&
o(r_n)
\eeas
as $n\to\infty$. 
\end{theorem}

\section{Some details of derivation of the expansion for the quadratic form}\label{241224-1}

In this section, we will give somewhat detailed proof of Theorem \ref{220829-3}, 
which was originally presented in \cite{Yoshida2010}. 

\subsection{Stochastic expansion} 
We will work with the It\^o stochastic integral equation (\ref{240221-1}). 
The following lemma gives a stochastic expansion of the targeted variable $Z_n$ of (\ref{240801-1}). 

\begin{en-text}
We write $f_t$ for $f(X_t)$, given function $f$. 
The It\^o decomposition of $\sigma_t=f(X_t)$ is 
denoted by 
\beas 
\sigma_t &=& f_0 + \int_0^t f^{[1]}_s dw_s
+\int_0^tf^{[0]}_sds. 
\eeas
This notation simplifies expressions though obviously 
$f^{[1]}_s=\partial_xf(X_s)\sigma(X_s)$ and 
\beas f^{[0]}_s=\partial_xf(X_s)[b(X_s)]+2^{-1}\partial_x^2f(X_s)[\sigma(X_s)^{\otimes2}]. \eeas 
\end{en-text}

\begin{lemma}\label{240803-2} $Z_n$ has the following stochastic expansion:
\beas 
Z_n &=& 
M^n_1+\frac{1}{\sqrt{n}}N_n, 
\eeas
where 
\beas 
M^n_t &=& 
\sqrt{n}\sum_{j=1}^n 2c_{t_{j-1}}\sigma_{t_{j-1}}^2
{\coloraka \int_{t_{j-1}\wedge t}^{t_j\wedge t}}
\int_{t_{j-1}}^sdw_rdw_s
\eeas
and 
\beas 
N_n &=& 
%
6n\sum_{j=1}^n c_{t_{j-1}}\sigma_{t_{j-1}}\sigma_{t_{j-1}}^{[1]}
\int_{t_{j-1}}^{t_j} \int_{t_{j-1}}^t \int_{t_{j-1}}^s
dw_udw_sdw_t
\\&&
+2\sum_{j=1}^n c_{t_{j-1}}b_{t_{j-1}}\sigma_{t_{j-1}}
\int_{t_{j-1}}^{t_j}dw_t
+
2n\sum_{j=1}^n c_{t_{j-1}}\sigma_{t_{j-1}} \sigma_{t_{j-1}}^{[1]} 
\int_{t_{j-1}}^{t_j}(t-t_{j-1}) dw_t
\\&&
{\coloraka
+
n^{-1}\sum_{j=1}^n  c_{t_{j-1}} b_{t_{j-1}}^2 
+n^{-1}\sum_{j=1}^n c_{t_{j-1}} \sigma_{t_{j-1}} b^{[1]}_{t_{j-1}} 
}
\\&&
{\coloro{-}}n\sum_{j=1}^n  c^{[1]}_{t_{j-1}} \sigma_{t_{j-1}}^2 
\int_{t_{j-1}}^{t_j} \int_{t_{j-1}}^t dw_sdt
\\&&
{\coloro{
-\frac{1}{2n}\sum_{j=1}^n c^{[0]}_{t_{j-1}}\sigma_{t_{j-1}}^2
-\frac{1}{n}\sum_{j=1}^n c^{[1]}_{t_{j-1}}\sigma_{t_{j-1}}\sigma^{[1]}_{t_{j-1}}
}}
+o_M(1). 
\eeas
Here $o_M(1)$ denotes a term of $o(1)$ as $n\to\infty$ 
with respect to $\bbD_{s,p}$-norms of any order. 
{\coloraka The families $\{M^n_t\}_{t\in[0,1],n\in\bbN}$ and $\{N_n\}_{n\in\bbN}$ are bounded in every $\bbD_{s,p}$-norm. }
\end{lemma}
It is possible to obtain the above lemma by somewhat long computations. 
See Section \ref{240803-1}.

\subsection{Random symbols}

\subsubsection{Adaptive random symbol}
\def\dla{\langle\!\langle}
\def\dra{\rangle\!\rangle}

The discrete filtration 
$\F^n=(\bar{\calf}^n_t)_{t\in  [0,1]}$ with 
$\bar{\calf}^n_t=\calf_{[nt]/n}$ will be used. 
The predictable quadratic covariation for $\F^n$ is denoted by 
$\dla\cdot,\cdot\dra$. Note that $\dla\cdot,\cdot\dra$ depends on $n$. 
Let 
$H_1(x)=x$ and $H_2(x)=(x^2-1)/\sqrt{2}$. 
Denote $\Delta_j w=w_{t_j}-w_{t_{j-1}}$, which depends on $n$ as well as $j$. 
The discrete version of $M^n$ is given by 
\beas 
\bar{M}^{2,n}_t 
&=&
\frac{1}{\sqrt{n}}\sum_{j:t_j\leq t}\sqrt{2}a(X_{t_{j-1}}) H_2(\sqrt{n}\Delta_j w). 
\eeas

For $\dotc\>\!\!^n_t:=\sqrt{n}\Big(C^n_t - C^\infty_t\Big)$, we have 
\beas 
\dotc\>\!\!^n_t
&=& 
\sum_{j:t_j\leq t}\int_{t_{j-1}}^{t_j}\>4n\sqrt{n}a(X_{t_{j-1}})^2
\Bigg\{\Big(\int_{t_{j-1}}^sdw_r\Big)^2-(s-t_{j-1})\Bigg\}ds
\\&&
-
2\sqrt{n}\sum_{j:t_j\leq t} \int_{t_{j-1}}^{t_j} \Big(a(X_s)^2-a(X_{t_{j-1}})^2\Big)ds
+O_p\Big(\frac{1}{\sqrt{n}}\Big)
\\&=& 
\sum_{j:t_j\leq t}\int_{t_{j-1}}^{t_j}\>4n\sqrt{n}a(X_{t_{j-1}})^2
\Bigg\{\Big(\int_{t_{j-1}}^sdw_r\Big)^2-(s-t_{j-1})\Bigg\}ds
\\&&
-
2\sqrt{n}\sum_{j:t_j\leq t} \int_{t_{j-1}}^{t_j}
2a(X_{t_{j-1}})a'(X_{t_{j-1}})(w_s-w_{t_{j-1}})ds
+O_p\Big(\frac{1}{\sqrt{n}}\Big) 
\\&=& 
\sum_{j:t_j\leq t}\int_{t_{j-1}}^{t_j}\>4n\sqrt{n}a(X_{t_{j-1}})^2
\Bigg\{\Big(\int_{t_{j-1}}^sdw_r\Big)^2-(s-t_{j-1})\Bigg\}ds
+O_p\Big(\frac{1}{\sqrt{n}}\Big). 
\eeas
Here 
the supremum of ``$O_p(n^{-1/2})$'' in $t\in[0,1]$ is of $O_p(n^{-1/2})$ with respect to $L^p$-norms. 
Therefore, the principal part of $\dotc\>\!\!^n$ is $\F^n$-martingale 
\beas 
\bar{M}^{\xi,n}_t 
&=&
\sum_{j:t_{j-1}\leq t}\int_{t_{j-1}}^{t_j} 4n\sqrt{n}a(X_{t_{j-1}})^2
\Bigg\{\Big(\int_{t_{j-1}}^sdw_r\Big)^2-(s-t_{j-1})\Bigg\}ds
\eeas
The discrete version of $w$ is denoted by $\bar{w}^n_t=w_{[nt]/n}$.

By a similar argument, we have 
\beas 
\sqrt{n}(F_n-F_\infty) &\to^p& 0=\dotf_\infty. 
\eeas

%
%

In the present situation, $\tilde{W}_\infty(z)=0$ and $\tilde{F}_\infty(z)=0$. 
We need to identify the limit $(M_\infty,\dotc_\infty,N_\infty)$ to write the adaptive random symbol. 
The ``martingale part'' of $N_n$ with respect to $\F^n$ is given by 
\beas 
\dot{N}^n_t
&=&
6n\sum_{j:t_j\leq t} c_{t_{j-1}}\sigma_{t_{j-1}}\sigma_{t_{j-1}}^{[1]}
\int_{t_{j-1}}^{t_j} \int_{t_{j-1}}^t \int_{t_{j-1}}^s
dw_udw_sdw_t
\\&&
+2\sum_{j:t_j\leq t}  c_{t_{j-1}}b_{t_{j-1}}\sigma_{t_{j-1}}
\int_{t_{j-1}}^{t_j}dw_t
+
2n\sum_{j:t_j\leq t}  c_{t_{j-1}}\sigma_{t_{j-1}} \sigma_{t_{j-1}}^{[1]} 
\int_{t_{j-1}}^{t_j}(t-t_{j-1}) dw_t
\\&&
{\coloraka
-n\sum_{j:t_j\leq t}   c^{[1]}_{t_{j-1}} \sigma_{t_{j-1}}^2 
\int_{t_{j-1}}^{t_j} \int_{t_{j-1}}^t dw_sdt. 
}
\eeas
%
%
Then
\beas 
\dla \bar{w}^n,\bar{w}^n \dra_t &=& \frac{[nt]}{n}\to t
\\
\dla \bar{M}^{2,n},\bar{M}^{2,n} \dra_t 
&=& 
\frac{2}{n}\sum_{j:t_{j-1}\leq t} a(X_{t_{j-1}})^2 
\to^p 
2\int_0^t a(X_s)^2 ds
\\
\dla \bar{M}^{\xi,n},\bar{M}^{\xi,n} \dra_t 
&=&
\frac{16}{3n}\sum_{j:t_{j-1}\leq t} a(X_{t_{j-1}})^4
\to^p 
\frac{16}{3}\int_0^t a(X_s)^4ds 
\\ 
{\coloraka \dla\bar{w}^n,\dot{N}\dra_t} &{\coloraka \to^p}&{\coloraka \int_0^t k_sds } 
\\
\dla \bar{w}^n,\bar{M}^{k,n} \dra_t &=& 0
\sskip(k=2,\xi)
\\
\dla \bar{M}^{2,n},\bar{M}^{\xi,n} \dra_t &=& 
\frac{{\coloro{8}}}{3n}\sum_{j:t_{j-1}\leq t} a(X_{t_{j-1}})^3
\to^p 
\frac{{\coloro{8}}}{3}\int_0^t a(X_s)^3 ds,
\\
\dla \bar{M}^{2,n},\dot{N}^n \dra_t &\to^p&0,
\\
\dla \bar{M}^{\xi,n},\dot{N}^n \dra_t &\to^p&0,
\\
\dla \dot{N}^n,\dot{N}^n \dra_t &\to^p& \int_0^t q_s^2 ds
\eeas
as $n\to\infty$ for each $t\in[0,1]$, 
where $\bbR_+$-valued process $q_t$ 
takes the form 
\beas 
q_t^2 &=& p(c_t,c^{[1]}_t,b_t,\sigma_t,\sigma^{[1]}_t)
\eeas
for some polynomial $p$; 
\begin{en-text}
is defined by 
\beas 
q_t^2 &=&
\frac{22}{3}c_t^2\sigma_t^2(\sigma^{[1]}_t)^2
+\frac{1}{3}(c^{[1]}_t)^2\sigma_t^4
+\frac{2}{3}c_tc^{[1]}_t\sigma_t^3\sigma^{[1]}_t
+\frac{28}{3}c_t^2b_t^2\sigma_t^2
+\frac{16}{3}c_t^2b_t\sigma_t^2\sigma^{[1]}_t
+\frac{10}{3}c_tb_t\sigma_t^3c^{[1]};
\eeas
\end{en-text}
it is possible to give an explicit expression of $p$, 
however we do not need the precise form of $q_t$ later. 
The orthogonality of $\bar{M}^{2,n}$, $\bar{M}^{\xi,n}$ and $\dot{N}^n$ 
to any bounded martingale orthogonal to $w$ is obvious, thus 
with the representation of $\dotc\>\!\!^n_t$, 
$M^n_t$ and $N_n$, we obtain
\beas 
(M_\infty,\dotc_\infty,N_\infty)
&=^d&
\Big( \int_0^1 \sqrt{2}a(X_s)dB_s,
\int_0^1 \frac{{\coloro{4}}\sqrt{2}}{3}a(X_s)^2dB_s
+\int_0^1 \frac{{\coloro{4}}}{3}a(X_s)^2dB'_s,
\\&&
{\coloraka \int_0^1 k_sdw_s+}\int_0^1 {\coloraka \sqrt{q_s^2-k_s^2}}dB''_s+\int_0^1 h_sds
\Big), 
\eeas
where $(B,B', B'')$ is a three-dimensional standard Wiener process, 
independent of $\calf$,  
defined on the extension $\bar{\Omega}$.  
Since 
\beas 
\tilde{N}_\infty(z) &=& {\coloraka \int_0^1 k_tdw_t+}\int_0^1 h_tdt, 
\eeas
the random symbol $\underline{\sigma}(z,iu,iv)$ is given by (\ref{240513-1}). 
Moreover we see Condition $[\bbA2]$ holds.

\subsubsection{Anticipative random symbol}

Let us find the anticipative random symbol $\overline{\sigma}(iu,iv)$. 
Recall that $\alpha(x)=a(x)^2$, 
\beas 
C_s^\infty=2\int_0^s\alpha(X_t)dt, \sskip
C_\infty=2\int_0^1 \alpha(X_t)dt, \sskip
F_\infty=\int_0^1\beta(X_t)dt\sskip\mbox{and}\sskip
W_\infty=0. 
\eeas 
{\coloraka To describe $\sigma_{s,s}(iu,iv)$ more precisely, we consider 
the random symbol $\sigma_{s,r}(iu,iv)$ that admits the expression}
\beas 
&&
\sigma_{s,r}(iu,iv)
\\&=&
u^2\int_r^s\alpha'(X_t)D_rX_tdt
\Big(- u^2\int_s^1\alpha'(X_t)D_sX_tdt+i\int_s^1\beta'(X_t)[v]D_sX_tdt \Big)
\\&&
+
\Big(-u^2\int_r^1\alpha'(X_t)D_rX_tdt+i\int_r^1\beta'(X_t)[v]D_rX_tdt \Big)
\Big(-u^2\int_s^1\alpha'(X_t)D_sX_tdt+i\int_s^1\beta'(X_t)[v]D_sX_tdt \Big)
\\&&
+
\Big(-u^2\int_s^1\{\alpha''(X_t)D_rX_tD_sX_t
+\alpha'(X_t)D_rD_sX_t\}dt
+i\int_s^1\{\beta''(X_t)[v]D_rX_tD_sX_t
+\beta'(X_t)[v]D_rD_sX_t\}dt
\Big)
\eeas
for $r\leq s$, 
where the prime $'$ stands for the derivative in $x_1\in\bbR$. 
The random symbol $\sigma_{s,s}$ is the limit $\lim_{r\up s}\sigma_{s,r}(iu,iv)$, that is, (\ref{26720828-2}), 
\begin{en-text}
\beas 
\sigma_{s,s}(iu,iv)
&=&
\Big(-u^2\int_s^1\alpha'(X_t)D_sX_tdt+i\int_s^1\beta'(X_t)[v]D_sX_tdt \Big)^2
\\&&
+
\Big(-u^2\int_s^1\{\alpha''(X_t)D_sX_tD_sX_t
+\alpha'(X_t)D_sD_sX_t\}dt
\\&&
+i\int_s^1\{\beta''(X_t)[v]D_sX_tD_sX_t
+\beta'(X_t)[v]D_sD_sX_t\}dt
\Big)
\eeas
\end{en-text}
and this gives the anticipative random symbol (\ref{26720828-1}).

\subsection{Nondegeneracy}

Here we will consider the nondegeneracy of $(M^n_t,F_\infty)$ in 
Malliavin's sense. 
{\coloroy{
Let 
\beas 
\eta_j(t) = \sqrt{n}(w(t_j\wedge t)-w(t_{j-1}\wedge t))
\eeas 
and 
\beas 
\xi_j(t) &=& 
n\bigg((w(t_j\wedge t)-w(t_{j-1}\wedge t))^2-(t_j\wedge t - t_{j-1}\wedge t)\bigg). 
\eeas
Then 
\beas 
M^n_t &=& 
\frac{1}{\sqrt{n}}\sum_{j=1}^n a(X\tjm)\xi_j(t).
\eeas
The representing density of the Malliavin derivative of $D^n_t$ is  
\beas 
D_rM^n_t &=& 
\sum_{j=1}^n
2{\coloroy{a({\coloraka X}_{t_{j-1}})}}
\eta_j(t) 1_{(t_{j-1}\wedge t,t_j\wedge t]}(r)
\\&&
+\frac{1}{\sqrt{n}}\sum_{j=1}^na'(X_{t_{j-1}})D_rX\tjm\xi_j(t)1_{\{r\leq t_{j-1}\leq t\}}
\\&=:&
D_1(n,t)_r+D_2(n,t)_r. 
\eeas
Now 
\beas 
D_2(n,t)_r 
&=& 
\frac{1}{\sqrt{n}}\sum_{j=2}^na'(X_{t_{j-1}})D_rX\tjm\xi_j(t)1_{\{r\leq t_{j-1}\leq t\}}
\\&=&
\frac{1}{\sqrt{n}}\sum_{\ell=1}^{n-1} a'(X_{t_{\ell}})\xi_{\ell+1}(t)D_rX_{t_{\ell}}1_{\{r\leq t_\ell\leq t\}}
\\&=&
\frac{1}{\sqrt{n}}\sum_{\ell=1}^{n-1} \sum_{j=1}^{\ell}a'(X_{t_{\ell}})\xi_{\ell+1}(t)D_rX_{t_{\ell}}1_{I_j(t)}(r)
\\&=&
\frac{1}{\sqrt{n}}\sum_{j=1}^{n-1}\bigg(\sum_{\ell=j}^{n-1} a'(X_{t_{\ell}})\xi_{\ell+1}(t)D_rX_{t_{\ell}}\bigg)1_{I_j(t)}(r)
\\&=&
\frac{1}{\sqrt{n}}\sum_{j=1}^{n-1}\bigg(\sum_{k=j+1}^n a'(X_{t_{k-1}})\xi_k(t)D_rX_{t_{k-1}}\bigg)1_{I_j(t)}(r),
\eeas
where $I_j(t)=(t_{j-1}\wedge t,t_j\wedge t]$.  
Hence the Malliavin covariance of $M^n_t$ is 
\beas 
\sigma_{11}(n,t):=\sigma_{M^n_t} 
&=&
\sum_{j=1}^n \int_{t_{j-1}\wedge t}^{t_j\wedge t}\bigg[2a({\coloraka X}_{t_{j-1}})\eta_j(t)
+\frac{1}{\sqrt{n}}\sum_{k=j+1}^na'(X_{t_{k-1}})D_rX_{t_{k-1}}\xi_k(t)\bigg]^2dr
\eeas
for $u\in\bbR$, 
where we read $\sum_{k=n+1}^n...=0$. 
Since 
\beas 
D_r F_\infty &=& \int_r^1\beta'(X_s)D_rX_sds, 
\eeas
we have 
\beas &&
\sigma_{12}(n,t)[v]
:=\langle M^n_t,F_\infty[v]\rangle_H
\\&=& 
\sum_{j=1}^n 
\int_{r=t_{j-1}\wedge t}^{t_j\wedge t}\int_{s=r}^1
\bigg[2a(X_{t_{j-1}})\eta_j(t)
+\frac{1}{\sqrt{n}}\sum_{k=j+1}^na'(X_{t_{k-1}})\xi_k(t)D_rX_{t_{k-1}}\bigg]
\beta'(X_s)D_rX_sdsdr[v]
\eeas
for $v\in\bbR^{d_1}$. Let 
\beas 
\sigma_{22}(t)[v^{\otimes2}]
&=&
\int_0^t\bigg[\int_r^1\beta'(X_s)D_rX_sds[v]\bigg]^2dr.
\eeas
and let 
\beas
\sigma(n,t)&=&
\l[
\begin{array}{cc} 
\sigma_{11}(n,t) & \sigma_{12}(n,t)^{{\coloraka \star}} \\ \sigma_{12}(n,t) & \sigma_{22}(t)
\end{array} 
\r]. 
\eeas
Let 
\beas 
\tilde{\sigma}_{11}(n,t)
&=&
\frac{1}{n}\sum_{j=1}^n \bigg[2a(X_{t_{j-1}})\eta_j(t)\bigg]^2
+\sum_{j=1}^n\int_{t_{j-1}\wedge t}^{t_j\wedge t} \bigg[\frac{1}{\sqrt{n}}\sum_{k=j+1}^na'(X_{t_{k-1}})\xi_k(t)D_rX_{t_{k-1}}\bigg]^2dr,
\eeas
\beas 
\tilde{\sigma}_{12}(n,t)[v]
&=& 
\sum_{j=1}^n \int_{r=t_{j-1}\wedge t}^{t_j\wedge t} \int_{s=r}^1 
\frac{1}{\sqrt{n}}\sum_{k=j+1}^na'(X_{t_{k-1}})\xi_k(t)D_rX_{t_{k-1}}\beta'(X_s)D_rX_sdsdr[v]
\eeas
and 
\beas 
\tilde{\sigma}(n,t)
=\l[
\begin{array}{cc} 
\tilde{\sigma}_{11}(n,t) & \tilde{\sigma}_{12}(n,t)^{{\coloraka \star}}  \\ \tilde{\sigma}_{12}(n,t) & \sigma_{22}(t)
\end{array} 
\r]. 
\eeas
We shall show 
\bea\label{220820-1} 
\big\| \sigma(n,t) - \tilde{\sigma}(n,t) \big\|_p &=& O(n^{-\half})
\eea
for every $p>1$ and $t\in[0,1]$ (in particular, for $t=1/2$). 
For this purpose, we need a lemma. 
Let $\cali$ denote the set of sequences $J^{(\nu)}=(J^{(\nu)}_{n,j})$ of 
multiple It\^o stochastic integrals  
taking the form
\beas 
J^{(\nu)}_{n,j} 
&=& n^{\frac{\nu}{2}}
\int_{t_{j-1}}^{t_j}dw_{s_1}a^{(\nu)}_{n,j,1}(s_1)
\int_{t_{j-1}}^{s_1}dw_{s_2}a^{(\nu)}_{n,j,2}(s_2)
\int_{t_{j-1}}^{s_2}
\cdots
\int_{t_{j-1}}^{s_\nu}dw_{s_\nu}a^{(\nu)}_{n,j,\nu}(s_\nu),
\eeas
where 
$\{a^{(\nu)}_{n,j,i};\>i=1,...,\nu$, $j=1,...,n$, $n\in\bbN\}$ is a family of progressively measurable processes,  
bounded in $\bbD^\infty=\cap_{s,p}\bbD_{s,p}$. 
In the following lemma, 
$J^{(\nu_1)}_{n,j_1}\cdots J^{(\nu_m)}_{n,j_m}$, 
$J^{(\mu_1)}_{n,k_1}\cdots J^{(\mu_m)}_{n,k_m}$ 
are in $\cali$, and each of them has $a^{(*)}_{n,j,i}$ which may possibly differ from 
those of other indices $\nu$'s and $\mu$'s even if 
the values of indices coincide each other.  
%
%
%
\begin{lemma}\label{220816-1}
Suppose that 
\beas 
&&
{\coloraka 
\sup_{j\in\{1,...,n\},\>n\in\bbN,\>\gamma\in\{0,1,...,m\},\atop\coloraka r_1,....,r_\gamma,s\in[0,1]}
\big\| |D_{r_1,....,r_\gamma}{\coloraka \calo}|\big\|_p
<\infty
}
\eeas
{\coloraka for all $\calo=a^{(*)}_{n,j}$, $b^{(*)}_{n,j}$ and $a^{(*)}_{n,j,i}(s)$, and } 
for every $p>1$.\footnote{$\gamma=0$ denotes the case with no derivative.} 
Then 
\bd
\item[(a)] 
Suppose that $b^{(d)}_{n,k}$ are $\calf_{t_{k-1}}$-measurable. 
Then for $\nu_1,...,\nu_m,\mu_1,...,\mu_q\in\bbN$, 
\beas 
\frac{1}{n^m}\sum_{j_1,....,j_m=1}^n E\bigg[
a^{(1)}_{n,j_1}J^{(\nu_1)}_{n,j_1}\cdots a^{(m)}_{n,j_m}J^{(\nu_m)}_{n,j_m}
\bigg(\frac{1}{\sqrt{n}}\sum_{k_1=j_1+1}^nb^{(1)}_{n,k_1}J^{(\mu_1)}_{n,k_1}\bigg)
\cdots
\bigg(\frac{1}{\sqrt{n}}\sum_{k_q=j_m+1}^nb^{(q)}_{n,k_q}J^{(\mu_q)}_{n,k_q}\bigg)
\bigg]
&=&
O\bigg(\frac{1}{n^{m/2}}\bigg). 
\eeas
\item[(b)] 
For $\nu_1,...,\nu_m\in\bbN$, 
\beas 
\frac{1}{n^m}\sum_{j_1,....,j_m=1}^n E\bigg[
a^{(1)}_{n,j_1}J^{(\nu_1)}_{n,j_1}\cdots a^{(m)}_{n,j_m}J^{(\nu_m)}_{n,j_m}\bigg]
&=&
O\bigg(\frac{1}{n^{m/2}}\bigg). 
\eeas
\ed
{\coloraka The constants in the above estimates depend only on the given supremums. } 
\end{lemma}
\begin{en-text}
\proof 
First we will show (a). 
We use the $L^2([0,T])$-orthogonality bewteen $1_{(t_{j-1},t_j]}$ and $1_{(t_{k-1},t_k]}$ for $j\not=k$. 
If the number of single $j_*$'s is $\alpha$, then the outside summation has 
at most $n^\alpha\times n^{\frac{m-\alpha}{2}}$ terms of such type. 
The $\alpha$ times IBP-formula for those single $j_*$'s deduces the order $n^{-m/2}$ 
if $k_1,...,k_q$ are different from any of $j_*$'s; otherwise, we also get $n^{-1/2}$ 
in each  IBP-formula. 
Note also that the derivative of $b$'s do not change the form of ``martingale'', besides 
it gives $n^{-1/2}$. 
After all, total order becomes 
\beas 
n^{-m}\times n^\alpha\times n^{\frac{m-\alpha}{2}}\times n^{-\half \alpha}
=n^{-\frac{m}{2}}. 
\eeas
In a similar way, we can obtain (b). 
\qed\\
\end{en-text}

By using smoothness of appearing functions, apply Lemma \ref{220816-1} (a) to 
\beas 
&&a^{(c)}_{n,j}=a(X\tjm),\sskip
J^{\nu_c}_{n,j}=\eta_j(t),
\\&&b^{(c)}_{n,k}=n\>a'(X_{t_{k-1}})\int_{t_{j-1}\wedge t}^{t_j\wedge t}D_rX_{t_{k-1}}dr\sskip(j<k)\sskip\mbox{and}\sskip
J^{\mu_c}_{n,k}=\xi_k(t)
\eeas
to obtain 
\beas 
\sup_{t\in[0,1]}\big\|\sigma_{11}(n,t)-\tilde{\sigma}_{11}(n,t)\big\|_p&=&O\l(\frac{1}{\sqrt{n}}\r)
\eeas
Furthermore, Lemma \ref{220816-1} (b) applied to 
\beas
a^{(c)}_{n,j}
=
\mbox{element of }
2a(X\tjm)\times 
n\int_{r=t_{j-1}\wedge t}^{t_j\wedge t}\int_{s=r}^1\beta'(X_s)D_rX_sdsdr{\coloraka [v]}
\sskip\mbox{and}\sskip J^{(\nu_c)}_{n,j}=\eta_j(t)
\eeas
yields 
\beas 
\sup_{t\in[0,1]}\big\|\sigma_{12}(n,t)-\tilde{\sigma}_{12}(n,t)\big\|_p
&=& O\l(\frac{1}{\sqrt{n}}\r)
\eeas
as $n\to\infty$ for every $p>1$. 
Consequently, we obtain (\ref{220820-1}).

Let $m_n=n^{-1}\sum_{j=1}^n\eta_j(1/2)^2$. 
For a positive number $c_1$, define $s_n$ by 
\beas 
s_n &=& \half\det\bigg[\tilde{\sigma}\bigg(n,\half\bigg)+\psi\l(\frac{m_n}{2c_1}\r)I_{1+d_1}\bigg],
\eeas
where $\psi:\bbR\to[0,1]$ is a smooth function such that $\psi(x)=1$ if $|x|\leq1/2$ and $\psi(x)=0$ if $|x|\geq1$. 
Then $s_n\geq 2^{-1}$ if $m_n\leq c_1$, and $s_n\geq2^{-1}\det\tilde{\sigma}(n,1/2)$ otherwise. 
Thus, it suffices to show 
\bea\label{240825-1}
\sup_n E\big[1_{\{m_n\geq c_1\}}\big(\det\tilde{\sigma}(n,1/2)\big)^{-p}\big]&<&\infty\sskip(p>1)
\eea
for the nondegeneracy 
\bea\label{240816-11}
\sup_nE[s_n^{-p}]&<&\infty
\eea
for every $p>1$. 
Following precisely, e.g., the proof of Lemma 2.3.1 of \cite{Nualart2006}, 
in order to obtain (\ref{240825-1}), 
it is sufficient to show that 
for every $p>1$, there exists a constant $C_p$ such that 
\bea\label{240825-2} 
\sup_{{\bf u}\in\bbR^{1+d_1}:|{\bf u}|=1}P\big[m_n\geq c_1,\>\tilde{\sigma}(n,1/2)[{\bf u}^{\otimes2}]\leq\ep\big]
&\leq& C_p\ep^p
\eea
for all $\ep\in(0,1)$ and all $n\in\bbN$. 
[The reasoning there is valid even for the measures $1_{\{m_n\geq c_1\}}dP$ in place of $P$.]  
Here we use $L^{\infty-}$-boundedness of $\{\tilde{\sigma}(n,1/2)\}_{n\in\bbN}$. 
Besides, for a while we shall assume the nondegeneracy condition: for some constant $C_p$,  
\bea\label{240825-5}
\sup_{v\in\bbR^{d_1}:|v|=1}P\big[\sigma_{22}(1/2)[v^{\otimes2}]\leq\ep\big]&\leq&C_p\ep^p
\eea
for all $\ep\in(0,1)$. 
Suppose that $|{\bf u}|=1$ for  ${\bf u}=(u,v)\in\bbR^{1+d_1}$. 
For simplicity, we write $\tilde{\sigma}_{ij}$ for $\tilde{\sigma}_{ij}(n,1/2)$ and 
$\sigma_{22}$ for $\sigma_{22}(1/2)$. Let $p>1$. 
\begin{en-text}
When $|u|<\ep^{1/8}$, 
\beas &&
P\big[\tilde{\sigma}(n,1/2)[{\bf u}^{\otimes2}]\leq\ep\big]
\\&=&
P\big[
\tilde{\sigma}_{11}u^2+2\tilde{\sigma}_{12}[uv]+\sigma_{22}[v^{\otimes2}]\leq\ep\big]
\\&\leq&
P\bigg[2|\tilde{\sigma}_{12}[v]|>\frac{\sigma_{22}[v^{\otimes2}]}{2\ep^{1/8}}\bigg]
+P\big[\sigma_{22}[v^{\otimes2}]\leq2\ep\big]
\\&\leq&
P\big[4|\tilde{\sigma}_{12}|>\ep^{-1/16}\big]
+P\big[
\sigma_{22}[v^{\otimes2}]<\ep^{1/16}\big]
+P\big[\sigma_{22}[v^{\otimes2}]\leq2\ep\big]
\\&\simleq&
\ep^p
\eeas
uniformly in $n\in\bbN$ by (\ref{240825-5}) and the inequality $(1-\ep^{1/4})^{1/2}<|v|\leq1$. 
\end{en-text}
When $|v|<\ep^{1/8}$, 
\bea\label{240826-5} &&
P\big[m_n\geq c_1,\>\tilde{\sigma}(n,1/2)[{\bf u}^{\otimes2}]\leq\ep\big]
\nn\\&=&
P\big[m_n\geq c_1,\>
\tilde{\sigma}_{11}u^2+2\tilde{\sigma}_{12}[uv]+\sigma_{22}[v^{\otimes2}]\leq\ep\big]
\nn\\&\leq&
P\bigg[{\coloraka m_n\geq c_1,\>}
2|\tilde{\sigma}_{12}u|>\frac{\tilde{\sigma}_{11}u^2}{2\ep^{1/8}}\bigg]
+P\big[m_n\geq c_1,\>\tilde{\sigma}_{11}u^2\leq2\ep\big]
\nn\\&\leq&
P\big[4|\tilde{\sigma}_{12}|>\ep^{-1/16}{\coloraka |u|}\big]
+P\big[m_n\geq c_1,\>
\tilde{\sigma}_{11}u^2<\ep^{1/16}\big]
+P\big[m_n\geq c_1,\>\tilde{\sigma}_{11}u^2\leq2\ep\big]
\nn\\&\simleq&
\ep^p
\eea
uniformly in $n\in\bbN$ and ${\bf u}$ satisfying $|{\bf u}|=1$ and $|v|<\ep^{1/8}$, 
where $c_2:=\inf_x|a(x)|>0$, since $(1-\ep^{1/4})^{1/2}<|u|\leq1$ and 
\beas 
\tilde{\sigma}_{11}u^2 
&\geq& 
2c_1c_2^2
\eeas
for any $\ep\in(0,2^{-4})$ on the event $\{m_n\geq c_1\}$. 

We will assume 
$|v|\geq\ep^{1/8}$. 
{From} 
\beas 
\tilde{\sigma}_{11}u^2+2\tilde{\sigma}_{12}[uv]+\sigma_{22}[v^{\otimes2}]
&\geq&
\tilde{\sigma}_{11}^{-1}\bigg\{\tilde{\sigma}_{11}\sigma_{22}[v^{\otimes2}]-\big(\tilde{\sigma}_{12}[v]\big)^2\bigg\},
\eeas
it follows that
\bea\label{240825-12} &&
P\big[m_n\geq c_1,\>\tilde{\sigma}(n,1/2)[{\bf u}^{\otimes2}]\leq\ep\big]
\nn\\&=&
P\big[m_n\geq c_1,\>
\tilde{\sigma}_{11}u^2+2\tilde{\sigma}_{12}[uv]+\sigma_{22}[v^{\otimes2}]\leq\ep\big]
\nn\\&\leq&
P\big[\tilde{\sigma}_{11}>\ep^{-1/4}\big]
+
P\big[m_n\geq c_1,\>
\tilde{\sigma}_{11}\sigma_{22}[v^{\otimes2}]-\big(\tilde{\sigma}_{12}[v]\big)^2\leq\ep^{3/4}\big].
\eea
On the other hand, we have 
\bea\label{220820-2} &&
\tilde{\sigma}_{11}\sigma_{22}[v^{\otimes2}]-\big(\tilde{\sigma}_{12}[v]\big)^2
\nn\\&=&
\frac{1}{n}\sum_{j=1}^n \bigg[2a(X_{t_{j-1}})\eta_j(t)\bigg]^2\int_0^t\bigg[\int_r^1\beta'(X_s)D_rX_sds[v]\bigg]^2dr
\nn\\&&
+\sum_{j=1}^n\int_{t_{j-1}\wedge t}^{t_j\wedge t} \bigg[\frac{1}{\sqrt{n}}\sum_{k=j+1}^na'(X_{t_{k-1}})\xi_k(t)D_rX_{t_{k-1}}\bigg]^2dr
\int_0^t\bigg[\int_r^1\beta'(X_s)D_rX_sds[v]\bigg]^2dr
\nn\\&&
-\bigg\{
\sum_{j=1}^n \int_{r=t_{j-1}\wedge t}^{t_j\wedge t} 
\bigg(\frac{1}{\sqrt{n}}\sum_{k=j+1}^na'(X_{t_{k-1}})\xi_k(t)D_rX_{t_{k-1}}\bigg)
\bigg(\int_{s=r}^1 \beta'(X_s)D_rX_sds[v]\bigg)dr
\bigg\}^2
\nn\\&\geq&
\frac{1}{n}\sum_{j=1}^n \bigg[2a(X_{t_{j-1}})\eta_j(t)\bigg]^2\int_0^t\bigg[\int_r^1\beta'(X_s)D_rX_sds[v]\bigg]^2dr
\nn
\eea
for $t=1/2$, 
where we used the Schwarz inequality. 
Hence, 
\bea\label{240825-11}
\tilde{\sigma}_{11}\sigma_{22}[v^{\otimes2}]-\big(\tilde{\sigma}_{12}[v]\big)^2
&\geq&
4c_1c_2^2\> \sigma_{22}(1/2)[v^{\otimes2}]
\eea
on $\{m_n\geq c_1\}$. 
Combining (\ref{240825-5}) with scaling $v\mapsto\ep^{-1/8}v$ and (\ref{240825-11}), we obtain 
\bea\label{240826-1}&&
\sup_{v\in\bbR^{d_1}:\ep^{1/8}\leq|v|\leq1}
P\big[m_n\geq c_1,\>
\tilde{\sigma}_{11}\sigma_{22}[v^{\otimes2}]-\big(\tilde{\sigma}_{12}[v]\big)^2\leq\ep^{3/4}\big]
\nn\\&\leq&
\sup_{v\in\bbR^{d_1}:\ep^{1/8}\leq|v|\leq1}
P\big[4c_1c_2^2\> \sigma_{22}(1/2)[v^{\otimes2}]\leq\ep^{3/4}\big]
\nn\\&\leq&
\sup_{v\in\bbR^{d_1}:|v|=1}P\big[\sigma_{22}(1/2)[v^{\otimes2}]\leq(4c_1c_2^2)^{-1}\ep^{1/2}\big]
\nn\\&\simleq&\ep^p
\eea
for every $p>1$. 
By connecting (\ref{240826-1}) to (\ref{240825-12}), we obtain 
\bea\label{240826-6}
\sup_{{\bf u}\in\bbR^{1+d_1}:\>|{\bf u}|=1,\>\ep^{1/8}\leq|v|\leq1}
P\big[m_n\geq c_1,\>\tilde{\sigma}(n,1/2)[{\bf u}^{\otimes2}]\leq\ep\big]
&\simleq&
\ep^p
\eea
Thus from (\ref{240826-5}) and (\ref{240826-6}), we obtain (\ref{240825-2}) and hence (\ref{240816-11})  
under the assumption (\ref{240825-5}).

We consider the $(1+d_1)$-dimensional process $X=(X_t)_{t\in[0,1]}$ satisfying the stochastic differential equation 
\beas 
d\bar{X}_t &=& V_0(\bar{X}_t)dt+V_1(\bar{X}_t)\circ dw_t,\sskip
\bar{X}_0=(X_0,0)
\eeas
in the Stratonovich form. 
Then the H\"ormander condition $[\bbH2]$ 
together with the compactness of $\mbox{supp}P^{X_0}$ 
ensures that for $t\in(0,1]$ and 
for every $p>1$, there exists a constant $C_p$ such that 
\bea\label{26720826-1}
\sup_{{\bf v}\in\bbR^{1+d_1}:|{\bf v}|=1}
P\bigg[\int_0^tD_s\bar{X}_1\otimes D_s\bar{X}_1ds[{\bf v}^{\otimes2}]\leq\ep\bigg]
&\leq& C_p\ep^p
\eea
for all $\ep\in(0,1)$. See Kusuoka and Stroock \cite{KusuokaStroock1984,KusuokaStroock1985}, 
 Ikeda and Watanabe \cite{IkedaWatanabe1989}, Nualart \cite{Nualart2006}. 
In particular, (\ref{240825-5}) follows from (\ref{26720826-1}) applied to ${\bf v}=(0,v)$ for $v\in\bbR^{d_1}$. 

We choose $c_1$ such that $2c_1<1/2=\lim_{n\to\infty}E[m_n]$. 
Since  
\beas
\sigma_{(M^n_t,F_\infty)}
=\l[
\begin{array}{cc} 
\sigma_{11}(n,t) & \sigma_{12}(n,t) \\ \sigma_{12}(n,t) & \sigma_{22}(1)
\end{array} 
\r]
\geq \sigma(n,t)
\eeas
by definition, we see that for every $K>0$, 
\bea\label{240826-15} 
\sup_{t\geq1/2}P\big[\det\sigma_{(M^n_t,F_\infty)}<s_n\big]
&\leq&
P\big[\det\sigma(n,1/2)<s_n\big]+O(n^{-K})
\nn\\&\leq&
P\big[\det\tilde{\sigma}(n,1/2)<2s_n\big]+O(n^{-K})
\nn\\&\leq&
P\big[m_n<2c_1\big]+O(n^{-K})
\nn\\&=&
O(n^{-K})
\eea
as $n\to\infty$. Here we used (\ref{220820-1}). 
Properties (\ref{240826-15}) and (\ref{240816-11}) verify $[\bbA3]$.

\subsection{Proof of Theorem \ref{220829-3}}

It is easy to verify Condition $[\bbA1]$ under $[\bbH1]$. 
Conditions $[\bbA2]$ and $[\bbA3]$ have been proved in the preceding sections. 
Now we can apply Theorem \ref{26720826-2} to obtain Theorem \ref{220829-3}.


\section{Proof of Lemma \ref{240803-2}}\label{240803-1}
\noindent 
\begin{en-text}
Notation. 

For functions $b^i:[0,1]\to\bbR$, 
\beas 
J_j(b^1(b^2(\cdots (b^k)\cdots)))_t
&=&
\int_{t\wedge t_{j-1}}^{t\wedge t_j} b^1_{s_1}\int_{t_{j-1}}^{s_2} b^2_{s_2}\cdots\int_{t_{j-1}}^{s_k}b^k_{s_k}dw_{s_k}\cdots dw_{s_2}dw_{s_1} 
\eeas
Write simiply 
\beas 
J(b^1(b^2(\cdots (b^k)\cdots)))_t
&=& 
J_j(b^1(b^2(\cdots (b^k)\cdots)))_t\sskip(t\in(t_{j-1},t_j])
\eeas

\beas 
I_j(b)_t&=&\int_{t\wedge t_{j-1}}^{t\wedge t_{j-1}}b_sds
\eeas
\beas 
I(b)_t&=&I_j(b)_t\sskip(t\in(t_{j-1},t_j])
\eeas
\end{en-text}

Let 
\beas 
&&f\one = \sigma,\sskip
f\two = b
\\&&
f^{(1,1)} = \sigma^{[1]}=\sigma'\sigma,\sskip
f^{(1,2)} = \sigma^{[0]}=\sigma'b+\half \sigma''\sigma^2.
\eeas
We have 
\bea\label{240813-1} &&
\sum_jc\tjm(X\tj-X\tjm)^2
\nn\\&=&
\sum_jc\tjm\intj(f\one_t)^2dt
+2\sum_jc\tjm\intj\l(\intjm^t f\one_sdw_s+\intjm^t f\two_sds \r)f\one_tdw_t
\nn\\&&
+2\sum_jc\tjm\intj\l(\intjm^t f\one_sdw_s+\intjm^t f\two_sds \r)f\two_tdt.
\eea
In the decomposition (\ref{240813-1}), 
\beas &&
\sum_jc\tjm\intj(f\one_t)^2dt
\\&=&
\frac{1}{2n}\sum_jc\tjm(f\one\tjm)^2
+\sum_jc\tjm\intj\intjm^t2f\one_s f^{(1,1)}_sdw_sdt
\\&&+\sum_jc\tjm\intj\intjm^t2f\one_s f^{(1,2)}_sdsdt
+\sum_jc\tjm\intj\intjm^t(f^{(1,1)}_s)^2dsdt
\\&=&
\frac{1}{2n}\sum_jc\tjm(f\one\tjm)^2
+\sum_j2c\tjm f\one\tjm f^{(1,1)}\tjm\intj\intjm^tdw_sdt
\\&&+\frac{1}{n^2}\sum_j c\tjm f\one\tjm f^{(1,2)}\tjm
+\frac{1}{2n^2}\sum_jc\tjm(f^{(1,1)}\tjm)^2
+o_p\l(\frac{1}{n}\r)
\eeas
as $n\to\infty$.
Next, 
\bea\label{240813-2} &&
2\sum_jc\tjm\intj\l(\intjm^t f\one_sdw_s+\intjm^t f\two_sds \r)f\one_tdw_t
\nn\\&=&
2\sum_jc\tjm f\one\tjm f\one\tjm\intj\intjm^t dw_sdw_t
+2\sum_jc\tjm f^{(1,1)}\tjm f\one\tjm\intj\intjm^t \int\tjm^udw_udw_sdw_t
\nn\\&&
+2\sum_jc\tjm f\two\tjm f\one\tjm\intj\intjm^t ds dw_t
+2\sum_jc\tjm f^{(1,1)}\tjm f\one\tjm\intj\l(\intjm^tdw_s \r)^2dw_t
+o_p\l(\frac{1}{n}\r)
\nn\\&=&
2\sum_jc\tjm f\one\tjm f\one\tjm\intj\intjm^t dw_sdw_t
+6\sum_jc\tjm f^{(1,1)}\tjm f\one\tjm\intj\intjm^t \int\tjm^udw_udw_sdw_t
\nn\\&&
+2\sum_jc\tjm f\two\tjm f\one\tjm\intj\intjm^t ds dw_t
+2\sum_jc\tjm f^{(1,1)}\tjm f\one\tjm\intj\intjm^tds dw_t
+o_p\l(\frac{1}{n}\r).
\eea
Further, 
\bea\label{240813-3} &&
2\sum_jc\tjm\intj\l(\intjm^t f\one_sdw_s+\intjm^t f\two_sds \r)f\two_tdt
\nn\\&=&
2\sum_jc\tjm f\one\tjm\intj\intjm^t dw_sf\two_tdt
+\frac{1}{n^2}\sum_jc\tjm (f\two\tjm)^2+o_p\l(\frac{1}{n}\r)
\nn\\&=&
2\sum_jc\tjm f\one\tjm f\two\tjm\intj\intjm^t dw_sdt
+2\sum_jc\tjm f\one\tjm f^{(2,1)}\tjm\intj\l(\intjm^t dw_s\r)^2dt
\nn\\&&
+\frac{1}{n^2}\sum_jc\tjm (f\two\tjm)^2+o_p\l(\frac{1}{n}\r)
\nn\\&=&
2\sum_jc\tjm f\one\tjm f\two\tjm\intj\intjm^t dw_sdt
+4\sum_jc\tjm f\one\tjm f^{(2,1)}\tjm\intj\intjm^t\intjm^sdw_u dw_sdt
\nn\\&&
+\frac{1}{n^2}\sum_jc\tjm f\one\tjm f^{(2,1)}
+\frac{1}{n^2}\sum_jc\tjm f\two\tjm f\two\tjm+o_p\l(\frac{1}{n}\r)
\nn\\&=&
2\sum_jc\tjm f\one\tjm f\two\tjm\intj\intjm^t dw_sdt
\nn\\&&
+\frac{1}{n^2}\sum_jc\tjm f\one\tjm f^{(2,1)}
+\frac{1}{n^2}\sum_jc\tjm (f\two\tjm)^2+o_p\l(\frac{1}{n}\r).
\eea

Let 
\beas 
c\one &=& c^{[1]}=c'\sigma=c'f\one,
\\
c\two &=& c^{[0]}=c'b+\half c''\sigma^2=c'f\two+\half c''(f\one)^2.
\eeas
Since
\beas 
\sigma_t &=& f\one(X_t)
\ =\ 
f\one\tjm +\intjm^tf^{(1,1)}_sdw_s+\intjm^tf^{(1,2)}_sds,
\eeas
we have 
\beas 
\sigma_t^2
&=& 
(f\one\tjm)^2+\intjm^t 2f\one_sf^{(1,1)}_sdw_s+\intjm^t 2f\one_sf^{(1,2)}_sds
+\intjm(f^{(1,1)}_s)^2ds.
\eeas
Now
\bea\label{240813-4} &&
\int_0^1 c_t\sigma^2_t dt 
=
\sum_j\intj c_t\sigma_t^2dt
\nn\\&=&
\sum_j\intj\bigg\{c\tjm+\intjm^tc\one_sdw_s+\intjm^tc\two_sds\bigg\}
\bigg\{(f\one\tjm)^2+\intjm^t2f\one_sf^{(1,1)}_sdw_s
\nn\\&&
+\intjm^t2f\one_sf^{(1,2)}_sds
+\intjm^t(f^{(1,1)}_s)^2ds\bigg\}dt
\nn\\&=&
\frac{1}{n}\sum_jc\tjm(f\one\tjm)^2
+\sum_j2c\tjm f\one\tjm f^{(1,1)}\tjm\intj\intjm^tdw_sdt
+\frac{1}{n^2}\sum_jc\tjm f\one\tjm f^{(1,2)}\tjm
+\frac{1}{2n^2}\sum_jc\tjm ( f^{(1,1)}\tjm)^2
\nn\\&&
+\sum_jc\one\tjm(f\one\tjm)^2\intj\intjm^tdw_sdt
\nn\\&&
+\sum_j\intj\intjm^tc\one_sdw_s \intjm^t2f\one_sf^{(1,1)}_sdw_s dt
+\frac{1}{2n^2}\sum_jc\two\tjm(f\one\tjm)^2
+o_p\l(\frac{1}{n}\r)
\nn\\&=&
\frac{1}{n}\sum_jc\tjm(f\one\tjm)^2
+\sum_j2c\tjm f\one\tjm f^{(1,1)}\tjm\intj\intjm^tdw_sdt
+\frac{1}{n^2}\sum_jc\tjm f\one\tjm f^{(1,2)}\tjm
+\frac{1}{2n^2}\sum_jc\tjm (f^{(1,1)}\tjm)^2
\nn\\&&
+\sum_jc\one\tjm(f\one\tjm)^2\intj\intjm^tdw_sdt
\nn\\&&
+\frac{1}{n^2}\sum_jc\one\tjm f\one\tjm f^{(1,1)}\tjm 
+\frac{1}{2n^2}\sum_jc\two\tjm(f\one\tjm)^2
+o_p\l(\frac{1}{n}\r).
\eea
{From} (\ref{240813-1})-(\ref{240813-4}) we obtain 
\beas &&
\sum_jc\tjm(X\tj-X\tjm)^2 - \int_0^1c_t\sigma_t^2dt
\\&=&
\frac{1}{2n}\sum_jc\tjm(f\one\tjm)^2
+\sum_j2c\tjm f\one\tjm f^{(1,1)}\tjm\intj\intjm^tdw_sdt
\\&&+\frac{1}{n^2}\sum_j c\tjm f\one\tjm f^{(1,2)}\tjm
+\frac{1}{2n^2}\sum_jc\tjm(f^{(1,1)}\tjm)^2
\\&&
+2\sum_jc\tjm f\one\tjm f\one\tjm\intj\intjm^t dw_sdw_t
+6\sum_jc\tjm f^{(1,1)}\tjm f\one\tjm\intj\intjm^t \int\tjm^udw_udw_sdw_t
\nn\\&&
+2\sum_jc\tjm f\two\tjm f\one\tjm\intj\intjm^t ds dw_t
+2\sum_jc\tjm f^{(1,1)}\tjm f\one\tjm\intj\intjm^tds dw_t
\\&&
+2\sum_jc\tjm f\one\tjm f\two\tjm\intj\intjm^t dw_sdt
\nn\\&&
+\frac{1}{n^2}\sum_jc\tjm f\one\tjm f^{(2,1)}
+\frac{1}{n^2}\sum_jc\tjm (f\two\tjm)^2
\\&&
-\bigg\{
\frac{1}{n}\sum_jc\tjm(f\one\tjm)^2
+\sum_j2c\tjm f\one\tjm f^{(1,1)}\tjm\intj\intjm^tdw_sdt
+\frac{1}{n^2}\sum_jc\tjm f\one\tjm f^{(1,2)}\tjm
+\frac{1}{2n^2}\sum_jc\tjm (f^{(1,1)}\tjm)^2
\nn\\&&
+\sum_jc\one\tjm(f\one\tjm)^2\intj\intjm^tdw_sdt
\nn\\&&
+\frac{1}{n^2}\sum_jc\one\tjm f\one\tjm f^{(1,1)}\tjm 
+\frac{1}{2n^2}\sum_jc\two\tjm(f\one\tjm)^2
\bigg\}
+o_p\l(\frac{1}{n}\r)
\\&=&
2\sum_jc\tjm f\one\tjm f\one\tjm\intj\intjm^t dw_sdw_t
\\&&
+\bigg\{
6\sum_jc\tjm f^{(1,1)}\tjm f\one\tjm\intj\intjm^t \int\tjm^udw_udw_sdw_t
\nn\\&&
+2\sum_jc\tjm f\two\tjm f\one\tjm\intj\intjm^t ds dw_t
+2\sum_jc\tjm f^{(1,1)}\tjm f\one\tjm\intj\intjm^tds dw_t
\\&&
+2\sum_jc\tjm f\one\tjm f\two\tjm\intj\intjm^t dw_sdt
+\frac{1}{n^2}\sum_jc\tjm ( f\two\tjm)^2
+\frac{1}{n^2}\sum_jc\tjm f\one\tjm f^{(2,1)}\tjm
\\&&
-\sum_jc\one\tjm(f\one\tjm)^2\intj\intjm^tdw_sdt
\nn\\&&
-\frac{1}{2n^2}\sum_jc\two\tjm(f\one\tjm)^2
-\frac{1}{n^2}\sum_jc\one\tjm f\one\tjm f^{(1,1)}\tjm 
\bigg\}
+o_p\l(\frac{1}{n}\r).
\eeas
Thus we obtained 
\beas &&
\sqrt{n}\bigg(\sum_jc\tjm(X\tj-X\tjm)^2 - \int_0^1c_t\sigma_t^2dt\bigg)
\\&=&
\sqrt{n}\sum_j2c\tjm\sigma\tjm^2\intj\intjm^t dw_sdw_t
\\&&
+\frac{1}{\sqrt{n}}\bigg\{
6n\sum_jc\tjm\sigma\tjm \sigma^{[1]}\tjm\intj\intjm^t \int\tjm^udw_udw_sdw_t
\nn\\&&
+2n\sum_jc\tjm b\tjm \sigma\tjm\intj\intjm^t ds dw_t
+2n\sum_jc\tjm \sigma\tjm \sigma^{[1]}\tjm\intj\intjm^tds dw_t
\\&&
+2n\sum_jc\tjm \sigma\tjm b\tjm\intj\intjm^t dw_sdt
+\frac{1}{n}\sum_jc\tjm b\tjm^2
+\frac{1}{n}\sum_jc\tjm \sigma\tjm b^{[1]}\tjm
\\&&
-n\sum_jc^{[1]}\tjm\sigma\tjm^2\intj\intjm^tdw_sdt
\nn\\&&
-\frac{1}{2n}\sum_jc^{[0]}\tjm\sigma\tjm^2
-\frac{1}{n}\sum_jc^{[1]}\tjm \sigma\tjm \sigma^{[1]}\tjm 
\bigg\}
+o_p\l(\frac{1}{n}\r)
\\&=&
\sqrt{n}\sum_j2c\tjm\sigma\tjm^2\intj\intjm^t dw_sdw_t
\\&&
+\frac{1}{\sqrt{n}}\bigg\{
6n\sum_jc\tjm\sigma\tjm \sigma^{[1]}\tjm\intj\intjm^t \int\tjm^udw_udw_sdw_t
\nn\\&&
+2\sum_jc\tjm b\tjm \sigma\tjm\intj dw_t
+2n\sum_jc\tjm \sigma\tjm \sigma^{[1]}\tjm\intj\intjm^tds dw_t
\\&&
+\frac{1}{n}\sum_jc\tjm b\tjm^2
+\frac{1}{n}\sum_jc\tjm \sigma\tjm b^{[1]}\tjm
\\&&
-n\sum_jc^{[1]}\tjm\sigma\tjm^2\intj\intjm^tdw_sdt
\nn\\&&
-\frac{1}{2n}\sum_jc^{[0]}\tjm\sigma\tjm^2
-\frac{1}{n}\sum_jc^{[1]}\tjm \sigma\tjm \sigma^{[1]}\tjm 
+o_p(1)
\bigg\}.
\eeas

Moreover, we can obtain $\bbD_{s,p}$-estimates of the residual term as well as each term on the right-hand side. 
This completes the proof of Lemma \ref{240803-2}. 




\begin{en-text}
{\bf $\bigg[$} 
The following lemma does not hold and we do not use it!! 
\begin{lemma}\label{220816-2}
\beas 
\frac{1}{n^m}\sum_{j_1,....,j_m}^n E\bigg[
a^{(1)}_{n,j_1}\cdots a^{(m)}_{n,j_m}
\bigg(\frac{1}{\sqrt{n}}\sum_{k_1=j_1+1}^nb^{(1)}_{n,k_1}J^{(\mu_1)}_{n,k_1}\bigg)
\cdots
\bigg(\frac{1}{\sqrt{n}}\sum_{k_q=j_m+1}^nb^{(m)}_{n,k_q}J^{(\mu_m)}_{n,k_m}\bigg)
\bigg]
&=&
O\bigg(\frac{1}{n^{m/2}}\bigg). 
\eeas
\end{lemma}
{\bf $\bigg]$}

\koko
}}

}\coloroy{
By elementary argument, 
\beas 
\bigg\|\|D_1(n,1)\|_H^2-4\int_0^1a(w_t)^2dt\bigg\|_p &=& O(n^{-\half})
\eeas
for every $p>1$. 
[The first term on the right-hand side is dominating. $\leftarrow$ not true!]
On an extension of the original stochastic basis, 
there exists a Wiener process $B=(B_t)_{t\in[0,1]}$ independent of $\calf_1$ such that 
\beas 
\bigg(w,n^{-\half}\sum_{j=2}^{[n\cdot]} a'(w_{t_{j-1}})\xi_j\bigg)
&\to^d&
\bigg(w,\int_0^\cdot a'(w_t)\sqrt{2}dB_t\bigg)
\eeas
in $C([0,1];\bbR^2)$ as $n\to\infty$; 
this is a consequence of the theory of convergence of stochastic integrals 
by Jakbowskii et al. \cite{jak} or Kurtz and Protter \cite{Kur-Pro90}. 
Therefore, 
\beas 
\bigg\| \|D_2(n,1)\|_H^2 - \int_0^1 
\big(\int_r^1 a'(w_t)\sqrt{2}dB_t\big)^2 dr \bigg\|_p 
\koko
\eeas

}}
\koko

\beas 
D_rF_n &=& 
\frac{4}{n}\sum_{j=1}^n a'(w_{t_{j-1}})1_{(0,t_{j-1}]}(r).
\eeas

\beas 
|\langle DM^n_t,F_n\rangle |
&\leq&
\langle D_2(n,t) \rangle^\half \langle DF_n\rangle^\half
\\&\to&
0
\eeas
[An estimate of the convergence can be obtained. ]
\end{en-text}


\bibliographystyle{spmpsci}      
\bibliography{bibtex-20101014-20120827}   

\end{document}